\documentclass[a4paper,11pt]{article}

\usepackage{amsmath,amssymb,amsthm,fullpage}
\usepackage{enumerate}  
\usepackage{pgfplots}
\usepackage{subcaption}
\usepackage{authblk}

\newcommand{\RR}{\mathbb R}
\newcommand{\EE}{\mathbb E}
\newcommand{\hier}[3][X]{#1^{(#2)}_{#3}}
\newcommand{\salj}{\mathcal{F}}

\newcommand{\eps}{\varepsilon}
\newcommand{\ie}{i.e.\ }
\newcommand*{\abs}[1]{\lvert#1\rvert}
\newcommand*{\prob}[1]{\mathbb{P}(#1)}
\newcommand*{\bin}[1]{\operatorname{Bin}(#1)}

\newtheorem{theorem}{Theorem}[section]
\newtheorem{lemma}[theorem]{Lemma}
\newtheorem{proposition}[theorem]{Proposition}
\newtheorem{problem}[theorem]{Problem}

\theoremstyle{definition}
\newtheorem{case}{Case}

\theoremstyle{remark}
\newtheorem{remark}[theorem]{Remark}

\title{Competing types in preferential attachment graphs with community structure\thanks{Supported
    by UK Research and Innovation Future Leaders Fellowship MR/S016325/1 (J.H.) and European Research Council grant no.\ 883810 (J.H.), and by the Heilbronn Institute for Mathematical Research.}\ 
    }

\author[1]{John Haslegrave}
\author[2]{Jonathan Jordan}
\author[2]{Mark Yarrow}
\affil[1]{Lancaster University, United Kingdom.
    \texttt{j.haslegrave@lancaster.ac.uk}}
\affil[2]{University of Sheffield, United Kingdom. \texttt{jonathan.jordan,mark.yarrow@sheffield.ac.uk}}

\begin{document}
\maketitle
\begin{abstract}{We extend the two-type preferential attachment model of Antunovi\'c, Mossel 
	and R\'acz, where each new vertex takes its type according to a defined rule based on the types of its neighbours, to incorporate community structure, and investigate whether the proportions of vertices of each type synchronise between communities. The behaviour depends both on the choice of community structure and on the type assignment rule.  
	
	For essentially all cases where the single community model has more than one possible 
	limit, communities may fail to synchronise for weakly interacting communities. Even when 
	the single community model almost surely converges to a deterministic limit, 
	synchronisation is not guaranteed. However, we give natural conditions on the assignment rule 
	and, for two communities, on the structure, either of which will imply synchronisation to 
	this limit, and each of which is essentially best possible.
	
	We also give an example where the proportions of types almost surely do not converge, which is impossible in the single community model.}
	
	\medskip
	\noindent\textbf{Keywords:} {preferential attachment; vertex types; community structure; coexistence; random graphs} 
\end{abstract}

\section{Introduction}

In \cite{AMR}, Antunovi\'{c}, Mossel and R\'{a}cz introduced a model for preferential attachment graphs where each vertex is of one of a number of types, which may, for example, be thought of as brand preferences.  Each vertex chooses its type based on the types of its neighbours when it joins the network.  In this paper, we concentrate on the setting with two types, which we refer to as ``red'' and ``blue'', and standard preferential attachment. Antunovi\'{c}, Mossel and R\'{a}cz \cite{AMR} showed that in this setting the proportion of red vertices converges to a limit. In the ``linear model'' this limit is distributed with full support on $[0,1]$, whereas in non-linear models there are finitely many possible limits corresponding to fixed points of a particular polynomial $R$, defined later, which depends on the type assignment rule.

The purpose of this paper is to extend the analysis of \cite{AMR} to preferential attachment models with community structure. The model we use is the model for geometric preferential attachment graphs in Jordan \cite{geopref9}, in the special case where the space $S$ the graphs are embedded in is a finite set; each element of $S$ can then be thought of as a community, as described by Hajek and Sankagiri in \cite{hajek2018}.  The proportions of edge ends of each type within the different communities can then be thought of as a network of reinforced (possibly negatively) stochastic processes with an interaction described by the community structure.

There has been some interest in whether reinforced processes with an underlying community structure show \emph{synchronisation} in the sense that all communities show the same limiting behaviour.  For example Dai Pra, Louis and Minelli \cite{dai2014} show synchronisation for P\'{o}lya urns with a complete graph network structure, with each urn converging to the same limit, and Aletti, Crimaldi and Ghigletti \cite{aletti2017} show synchronisation for a more general reinforcement scheme based on a weighted network with an irreducible and diagonalisable matrix of weights.  Crimaldi, Louis and Minelli \cite{crimaldi2023} extend this to non-linear reinforcement, and as well as giving conditions for synchronisation also give examples with non-synchronisation.

In our setting it is thus a natural question whether different communities synchronise, or whether they can show different limiting behaviour for the proportions of types within them.  We will answer this question by showing (Theorem \ref{different}) that non-synchronisation is possible: where the polynomial $R$ defined by the type assignment rule has multiple fixed points, assuming some natural conditions on the type assignment rule, there is positive probability of different limits in different communities if the interaction between the communities is sufficiently weak.

We will also show (Theorems \ref{same1}, \ref{same2} and \ref{determinant}) that, again under some natural conditions on the interaction between the communities and the type assignment rule, this does not occur if the polynomial $R$ has only one fixed point (corresponding to only one limit having positive probability in the model of \cite{AMR}) and that in this the proportions of red in each community will all synchronise: they will converge to the same limit.  However, we will also show some examples, with perhaps less natural choices of community interaction and type assignment rule, where non-synchronisation is possible even in this case.

Based on the results in \cite{AMR} for a single community, it might be expected that for any fixed point of the polynomial $R$ there would be positive probability that the proportion of red vertices in each community converges to this fixed point.  It turns out that in the multi community set-up this is not necessarily the case, and indeed we will show an example (Proposition \ref{threecycle}) where the proportions of types in the different communities almost surely do not converge to limits.  However, we will give a result (Theorem \ref{stable}) with an extra condition on $R$ in the neighbourhood of the fixed point which ensures that there is positive probability of all communities converging to it; we also show that if there exists $z^*$ such that there is positive probability of all communities' proportions of red converging to $z^*$, then $z^*$ is a fixed point of $R$.

Finally, we will show (Theorem \ref{linear}) that, as long as there are enough non-zero interactions between communities, in the linear model (where $R(z)\equiv z$) synchronisation does occur, with all communities converging to the same random limit.

We mention some related work and some possible extensions.  The broadcasting problem studied by Addagio-Berry, Devroye, Lugosi and Velona \cite{broadcasting} on the preferential attachment tree is a special case of the model of \cite{AMR}, while Backhausz and Rozner \cite{BR19} consider a model where it is the edges rather than the vertices which have types.  Some further related work and open problems inspired by our model are discussed in section \ref{sec:open}.

The structure of the paper is as follows.  We define our framework in section \ref{sect:def} and state our results in section \ref{sect:results}, with proofs of the results in section \ref{sect:proofs}.  Finally in section \ref{examples} we consider some examples which illustrate our results, show some simulations, and give some more detailed calculations which are specific to those examples.

\section{Definitions and notation}\label{sect:def}
\subsection{Graph model}
Our graph model is that of \cite{geopref9} with the space $S$ considered to be a finite set of communities. It depends on three parameters: a positive integer $m$, a distribution $\mu$ on $S$, and a non-negative ``attractiveness'' matrix $A=(\alpha_{x,y})_{x,y\in S}$. Starting from a fixed initial graph, vertices $v_1,v_2,\ldots$ are added sequentially. Each vertex $v_t$ is assigned a community $c(v_t)\in S$, with the sequence $c(v_t)$ being i.i.d.\ of distribution $\mu$. When $v_t$ joins the graph, it forms $m$ edges to existing vertices.

The preferential attachment graph evolves from an initial graph $G_0$, with $n_0$ vertices, such that $G_0$ contains at least one vertex in each community and all vertices have positive degree. The graph $G_{t+1}$ is produced from $G_t$ by adding the vertex $v_{t+1}$ and 
$m$ edges from $v_{t+1}$ to $m$ vertices chosen independently at random (with replacement), with the probability of choosing a vertex $v_s$ being given by
\begin{equation}\label{attachment}\frac{\alpha_{c(v_s),c(v_{t+1})}\deg_{G_t}(v_s)}{\sum_{u\in V(G_t)}\alpha_{c(u),c(v_{t+1})}\deg_{G_t}(u)}.
\end{equation}
Note that, since the $m$ choices are independent, this may produce multiple edges between the same pair of vertices.  For notational convenience we assume the initial graph $G_0$ has $mn_0$ edges.

Thus the entry $\alpha_{x,y}$ is a multiplicative factor representing the relative attractiveness of an existing vertex in community $x$ to a new vertex in community $y$. (In many cases, $A$ will be a symmetric matrix, but this is not necessary in \cite{geopref9}.)  We allow an individual attractiveness coefficient $\alpha_{x,y}$ to be $0$, in which case a new vertex in community $y$ will never attach to a vertex in community $x$ (though if $\alpha_{y,x}>0$ edges may form between the two communities in the opposite way). However, for each community $y$ we require $\sum_x \alpha_{x,y}>0$, since a new vertex in community $y$ must attach somewhere. 

\subsection{Type assignment}
After joining the graph and forming edges, each vertex chooses one of two types, which we refer to as ``red'' and ``blue''. The choice of type is based on the types of its neighbours according to some rule, which may be random. This formulation allows the same type assignment rules as in \cite{AMR} to be considered, and we will use the same notation for the type assignment mechanism; in particular in the case with two types the probability that a new vertex becomes red, conditional on having $k$ red neighbours, is $p_k$.  Because the new vertex makes $m$ connections in total, the probability that a new vertex with $k$ blue neighbours becomes blue is $1-p_{m-k}$. Type assignment mechanisms with $p_k=1-p_{m-k}$ thus have symmetry between the two types.

The special case where $p_k=\frac{k}{m}$ for all $k$ is referred to as the \emph{linear model} in \cite{AMR}. The model studied in \cite{AMR} is as described above but with the graph growing according to a standard preferential attachment model, specifically the \emph{independent model} described in \cite{bergerpa}, which is the special case of the model of \cite{geopref9} with only one community;  we will refer to it in this paper as the \emph{single community} model.

We will associate the rule given by $p_0,\cdots,p_m$ with the polynomial given by 
\[R(z)=\sum_{i=0}^m\binom miz^i(1-z)^{m-i}p_i,\]
that is, $R(z)$ is the probability of a vertex becoming red if its neighbours are independently red with probability $z$. 

This polynomial allows us to describe the behaviour of the single community model. Except for the linear model where $R(z)\equiv z$, the possible limits for the proportion of red vertices are the fixed points of $R$; in the linear model the limit is a random variable with full support on $(0,1)$. This was proved in \cite{AMR}.

To understand the possible limits for the single community model further, for a fixed point $z\in (0,1)$ of $R$ we classify it as \emph{stable} if there exists $\eps>0$ such that $R(y)>y$ for $y\in (z-\eps,z)$ and $R(y)<y$ for $y\in (z,z+\eps)$, as \emph{unstable} if there exists $\eps>0$ such that $R(y)<y$ for $y\in (z-\eps,z)$ and $R(y)>y$ for $y\in (z,z+\eps)$, and as a \emph{touchpoint} if there exists $\eps>0$ such that $R(y)-y$ is either strictly positive or strictly negative on $(z-\eps,z)\cup(z,z+\eps)$.  Similarly $z=0$ is a stable fixed point of $R$ if $R(0)=0$ and there exists $\eps>0$ such that $R(y)<y$ for $y\in (0,\eps)$, and is an unstable fixed point of $R$ if $R(0)=0$ and there exists $\eps>0$ such that $R(y)>y$ if $y\in (0,\eps)$; similar definitions apply to a fixed point at $z=1$. We say that a stable fixed point $z$ is \emph{linearly stable} if the derivative $R'(z)<1$, and that an unstable fixed point $z$ is \emph{linearly unstable} if $R'(z)>1$. It is shown in \cite{AMR} that in the single community model the proportion of red edge-ends converges almost surely and that the limit is either a stable fixed point or a touchpoint of $R$, with all such points having positive probability of being the limit.

\section{Results and open problems}\label{sect:results}

Our results are concerned with the limiting proportion of red edge-ends in each community (if the limit exists). Equivalently this is the probability that a vertex
selected according to \eqref{attachment} is red, conditional on it being in that community, and as a result almost sure convergence of the proportions of edge ends in each community which are red implies the same for proportions of vertices in each community.

Write $\hier[Z]ti$ for the proportion of edge-ends in community $i$ at time $t$ which are at red vertices, that is 
\[\hier[Z]ti=\frac{\sum_{v\text{ red, }c(v)=i}\deg_{G_t}(v)}{\sum_{v:c(v)=i}\deg_{G_t}(v)}.
\]

Our first result shows that if the type assignment rule is such that $R$ has multiple linearly stable fixed points then we do get positive probability of different limits in different communities.  It will be useful to consider the following family of attractiveness matrices: given a matrix $A_1$, with all diagonal elements positive, let $A_0$ be a diagonal matrix whose diagonal elements are the same as those of $A_1$, and then for $\theta>0$ define $A_{\theta}=(1-\theta)A_0+\theta A_1$; changing the value of $\theta$ allows tuning of the relative strength of the inter-community and intra-community interactions.  Then the following theorem shows that for sufficiently small but positive $\theta$ there is positive probability that the proportions of red edge ends in each community converge to limits which are not the same in every community.  Furthermore, for any assignment of the linearly stable fixed points of $R$ to the different communities, there is positive probability that the limit in each community is close to the limit assigned to that community.

\begin{theorem}\label{different}Assume that the polynomial $R$ has at least two distinct linearly stable fixed points. For each $i\in N$, choose $z_i$ to be one of these points. Then for any choice of interaction matrix $A_1$ with all diagonal entries positive, and any $\eps>0$, there exists $\theta_\eps>0$ such that for $\theta<\theta_\eps$ the system with interaction matrix $A_{\theta}$ has positive probability that there exist limits $\tilde z_1,\ldots, \tilde z_N$, with $\abs{z_i-\tilde z_i}<\eps$, such that $\hier[Z]ti\to \tilde z_i$ as $t\to\infty$ for each $i$.

In particular, there is some $\theta_{\mathrm{crit}}>0$ such that for $\theta<\theta_{\mathrm{crit}}$ there is positive probability of convergence to different limits.\end{theorem} 
\begin{remark}\label{rem:weak-different}We also obtain a weaker form of this result for any $R$ with at least two distinct stable (\ie not necessarily linearly stable) fixed points, namely, that $\abs{\hier[Z]ti-z_i}<\eps$ for $t$ sufficiently large.
\end{remark}
Theorem \ref{different} cannot in general be extended to the case where $R$ has one stable fixed point and one touchpoint, even though in this case the single community model has more than one possible limit.  
\begin{proposition}\label{touchpoint}
	Suppose that the polynomial $R$ has a unique stable fixed point $z^*$, and some touchpoints. Let $A_1$ be any interaction matrix with all entries positive. Then for $\theta$ sufficiently small, in the system with interaction matrix $A_{\theta}$, convergence to different limits has probability zero.
\end{proposition}
Examples of this type where convergence to different limits has probability zero and others where convergence to different limits has positive probability will be considered in Section \ref{examples}.

Our second main result gives an almost complete classification of those type assignment rules that necessarily lead to the communities all having
the same limit almost surely, independently of the community structure. We divide this into two theorems, one showing sufficiency of our conditions and the other showing that marginally weaker conditions are also necessary.

\begin{theorem}\label{same1}
Suppose that the type assignment rule satisfies the following properties.
\begin{enumerate}
\item $R$ has a unique fixed point $z^*\in[0,1]$, except possibly for linearly unstable fixed points at $0$ and/or $1$
\item The only fixed points of $R(R(z))$ in $[0,1]$ are the fixed points of $R(z)$.
\end{enumerate}
Then for all communities $i$ we have $\hier[Z]ti \to z^*$ as $t\to\infty$, almost surely.
\end{theorem}

\begin{remark}We note that these conditions are satisfied by any increasing rule for which $R$ has a unique fixed point; here we say that a rule is \emph{increasing} if $R(z)$ is increasing in $z$. Since $R(z)$ may be expressed as $\EE(p_K)$ where $K\sim\bin{m,z}$, any rule for which $p_k$ is increasing in $k$ is necessarily an increasing rule; however, the rule being increasing is a strictly weaker condition. For example, when $m=3$ the rule given by $p_0=0, p_1=1,p_2=0,p_3=1$ is increasing, since $R(z)=z^3+3z(1-z)^2$, and so $R'(z)=3(1-2z)^2\geq0$.  All rules considered in \cite{AMR} in fact have $p_k$ increasing in $k$ and thus are increasing rules.
\end{remark}

In view of Theorem \ref{different}, condition 1 in Theorem \ref{same1} cannot be relaxed any further than allowing touchpoints or touchpoint-like unstable fixed points at $0$ and/or $1$; rules with such fixed points correspond to a set of measure zero in the rule space. Similarly we can show that condition 2 can only be relaxed to a similar extent. If $R$ satisfies condition 1 but not condition 2, then the smallest fixed point of $R(R(z))$ is either at $0$ or is stable (since $R(R(0))\geq 0$). Thus, in all but a measure-zero subset of cases, $R(R(z))$ has a linearly stable fixed point other than $z^*$. We show that Theorem \ref{same1} fails in all such cases.
\begin{theorem}\label{same2}
Assume that $z^*$ is the only fixed point of $R$ in $[0,1]$, except possibly for linearly unstable fixed points at $0$ and/or $1$. 
Suppose that the type assignment rule is such that $R(R(z))$ has a linearly stable fixed point $\tilde z\neq z^*$.
		Then for the community structure (with two communities) given by 
		\[A=\begin{pmatrix}0 & 1 \\ 1 & 0\end{pmatrix}\]
		and any $\mu_1,\mu_2>0$ with $\mu_1+\mu_2=1$, there is positive probability that there exist $z_1\neq z_2$ such that $\hier[Z]t1\to z_1$ and $\hier[Z]t2\to z_2$.
\end{theorem}

In the case where there are two communities, we can also give a necessary and sufficient condition on the matrix such that convergence to $z^*$ happens in both communities for any type assignment rule with a single fixed point of $R$. 

\begin{theorem}\label{determinant}
	Consider the two-community case, \ie $S=\{1,2\}$.
	\begin{enumerate}
		\item Suppose $\det A=\alpha_{1,1}\alpha_{2,2}-\alpha_{1,2}\alpha_{2,1} \geq 0$, and the type assignment rule is such that $z^*$ is the only fixed point of $R$ in $[0,1]$, except possibly for linearly unstable fixed points at $0$ and/or $1$.  Then for each $i\in\{1,2\}$ we have $\hier[Z]ti\to z^*$ as $t\to\infty$, almost surely.
		\item Conversely, suppose $\det A=\alpha_{1,1}\alpha_{2,2}-\alpha_{1,2}\alpha_{2,1} < 0$.  Then there exists a type assignment rule such that $R$ has a unique fixed point $z^*\in[0,1]$ but for which there is positive probability that there exist $z_1\neq z_2$ such that $\hier[Z]t1\to z_1$ and $\hier[Z]t2\to z_2$.
	\end{enumerate}
\end{theorem}

The criteria on $\det A=\alpha_{1,1}\alpha_{2,2}-\alpha_{1,2}\alpha_{2,1}$ here can be seen as measuring the tendency for vertices to connect to their own community ($\alpha_{1,1}\alpha_{2,2}-\alpha_{1,2}\alpha_{2,1}>0$, leading to an assortative community structure) or to connect outside their own community ($\alpha_{1,1}\alpha_{2,2}-\alpha_{1,2}\alpha_{2,1}<0$, leading to a disassortative community structure).

The methods used to prove Theorem \ref{same1} also tell us that for stable fixed points of $R$ where $R$ is locally increasing there is positive probability that all communities converge to them; in particular this will apply to any stable fixed point of $R$ when the rule is increasing.  Furthermore, we can show, without any conditions on $R$, that if there is a positive probability of synchronisation with the proportion of red in each community converging to $z^*$ then $z^*$ is a fixed point of $R$.  We combine these two facts into the following theorem.

\begin{theorem}\label{stable}\begin{enumerate}\item Let $z^{*}$ be a stable fixed point of the polynomial $R$ such that $R$ is increasing in some neighbourhood of $z^{*}$.  Then there is positive probability that $\hier[Z]ti\to z^*$ for all $i$.
\item Assume there is positive probability that $\hier[Z]ti\to z^*$ for all $i$.  Then $z^*$ is a fixed point of $R$.
\end{enumerate}
\end{theorem}

\begin{figure}[ht]
	\centering
	\scalebox{0.5}{\includegraphics{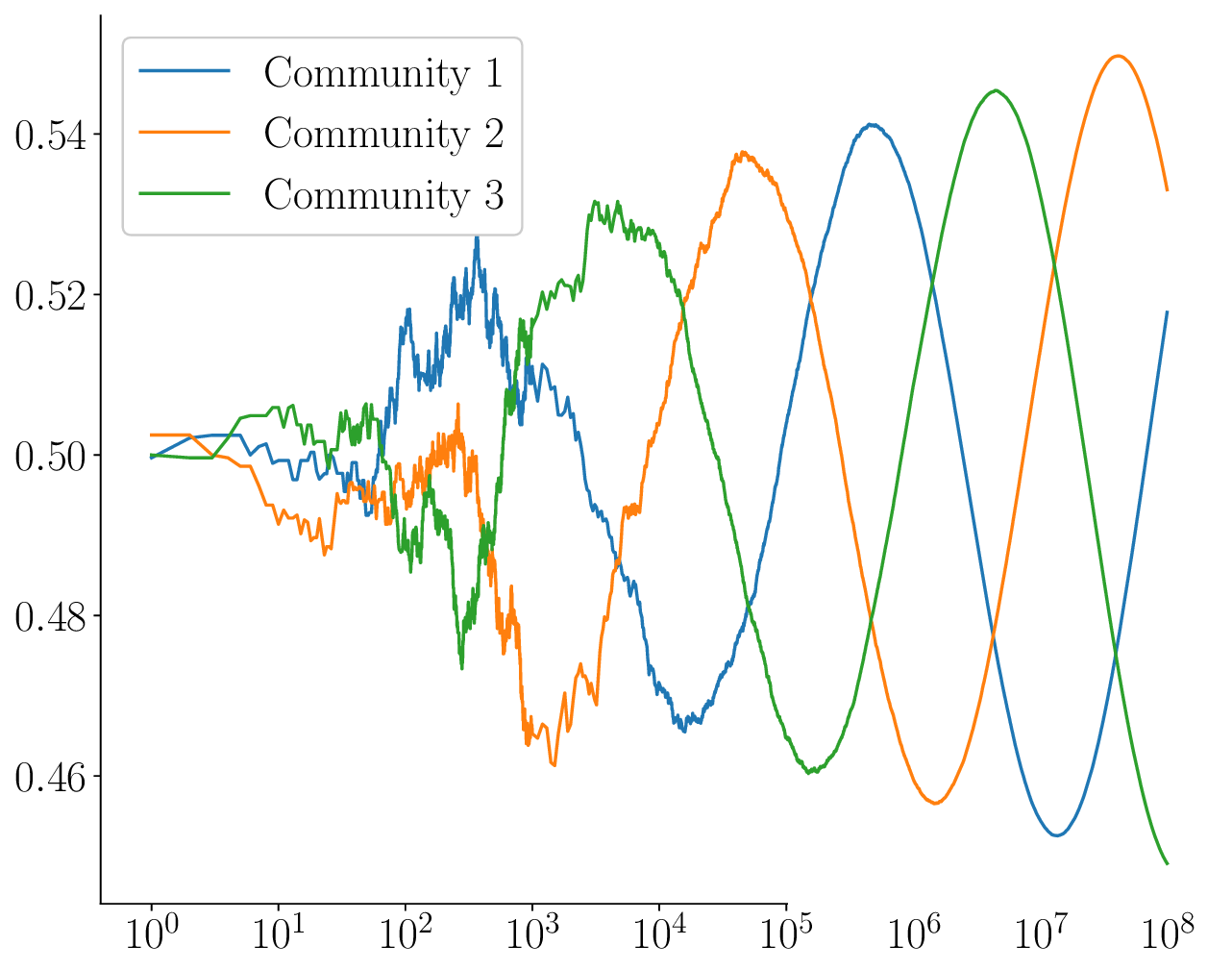}}
	\caption{The values of $\hier[Z]t1,\hier[Z]t2,\hier[Z]t3$ in the framework of Proposition \ref{threecycle} with $m=7$ plotted over time (log scale), showing apparent non-convergence.\label{fig:cycling}}
\end{figure}

Some condition beyond stability is necessary in the first part of Theorem \ref{stable}.  Indeed, it is possible to find examples which do not even converge to a limit; the following result gives such an example.
\begin{proposition}\label{threecycle}
	Let $S=\{1,2,3\}$ and \[A=\begin{pmatrix} 0 & 1 & 0 \\ 0 & 0 & 1 \\ 1 & 0 & 0\end{pmatrix}\]
	(representing three communities in a cycle, with each only influenced by the one clockwise of it), with $\mu_i=\frac13$ for $i=1,2,3$.  Let $m$ be an odd integer, and let the $p_k$ be given by the minority rule: $p_k=1$ for $0\leq k \leq \frac{m-1}{2}$ and $p_k=0$ for $\frac{m+1}{2} \leq k \leq m$.  Then for $m\geq 7$ almost surely $\bigl(\hier[Z]t1,\hier[Z]t2,\hier[Z]t3\bigr)$ does not converge to a limit. \end{proposition}
Figure \ref{fig:cycling} shows a simulation of proportions of red in each community over time in the framework of Proposition \ref{threecycle} for $m=7$.  It appears that the process is converging instead to a limit cycle.

Finally, we consider the linear model.  Here we have the following result, which shows that in this setting we do not get different limits in different communities for typical community structures.  To state it we need to consider the directed graph $\Gamma$ whose vertex set is $S$ and which has an edge from $i$ to $j$ if and only if $\alpha_{j,i}>0$.
\begin{theorem}\label{linear}
	Assume that $\Gamma$ is such that for any two vertices $i$ and $j$ there is a vertex $h$ such that there are directed paths from both $i$ and $j$ to $h$.  (We allow $h=i$ or $h=j$, and treat any vertex $i$ as having a directed path to itself.)  Then in the linear model with $p_k=\frac{k}{m}$ for all $k$, there exists a random limit $M$ such that $\hier[Z]ti\to M$ for all $i$.
\end{theorem}

\begin{remark}The condition on $\Gamma$ is equivalent to stating that a random walk on $\Gamma$ has a unique stationary distribution.  It is also equivalent to $\Gamma$ with every edge reversed being \emph{quasi strongly connected} as defined, for example, in \cite[Chapter 5]{MR1157371}.\end{remark}

\subsection{Open problems}\label{sec:open}

We suggest here a number of open problems and extensions which are left for further study.

\begin{problem}What happens if there are more than two types?\end{problem}
It is of course possible to define a version of this model with more than two types, and this was considered for the single community case in \cite[Section 3]{AMR}, although that paper focussed mainly on the two type case.  With more than two types there is potential for more complicated behaviour, and Haslegrave and Jordan \cite{threetypes} gave an example with non-convergence in that case.  We might therefore expect that there are interesting examples with multiple communities and more than two types, but we leave this for future work and focus on two types in this paper.

\begin{problem}
	The example in Proposition \ref{threecycle} shows that non-convergence is possible with two types, but with more than two communities. Is there an example with only two communities and two types which fails to converge with positive probability?
\end{problem}

\begin{problem}When there are two communities, Theorem \ref{determinant} gives a necessary and sufficient condition on the community structure for synchronisation to occur for any type assignment rule for which $R$ has a single fixed point.  What is the corresponding condition for more than two communities?\end{problem}
	
It is clear that the answer cannot be a simple criterion that $\det A\geq 0$, as the example \[A=\begin{pmatrix} 0 & 1 & 0 & 0 \\ 1 & 0 & 0 & 0 \\ 0 & 0 & 0 & 1 \\ 0 & 0 & 1 & 0\end{pmatrix}\] has positive determinant but would produce two essentially independent two-community structures, both of which have negative determinant.  A possible conjecture is that the criterion should be that all eigenvalues of $A$ have non-negative real part.

\begin{problem}Is non-synchronisation possible for continuous structure?\end{problem}
The model of \cite{geopref9} is more generally a geometric preferential attachment model with vertices embedded in a metric space $S$ as in Flaxman, Frieze and Vera \cite{FFV1,FFV2} and Manna and Sen \cite{manna}.  Given our results for the case where $S$ is finite, a natural question is whether different limits may be possible in different parts of $S$ when $S$ is an uncountable metric space such as a torus.  A degenerate case of geometric preferential attachment is the online nearest neighbour graph, as discussed in \cite{jordanwade} and \cite{manna} (where it is considered as the ``$\alpha=-\infty$'' case), where new vertices simply connect to the nearest vertex already present in the graph.   A simple type assignment process on the online nearest neighbour graph, where each vertex takes the type of its nearest neighbour when it joins the graph, is equivalent to the partitions model studied by Aldous in \cite{aldouspartitions} and by Basdevant, Blanc, Curien and Singh in the recent preprint \cite{basdevant2023}, where it is shown that with two initial seed vertices of different colours the process on $[0,1]^d$ converges in the Hausdorff sense to a random partition of the underlying set with a frontier which has a Hausdorff dimension strictly between $d-1$ and $d$; this might suggest that similar phenomena might occur for at least some geometric preferential attachment models, especially those whose behaviour is shown to be close to that of the online nearest neighbour graph in \cite{jordanwade}.

\section{Proofs}\label{sect:proofs}
\subsection{Stochastic approximation}
The methods used in \cite{AMR} rely heavily on the theory of stochastic approximation, relating the dynamics of the proportion of edge ends which are red to those of the differential equation driven by the polynomial $P(z)=(R(z)-z)/2$.  In this section we show how to set the process up as a multi-dimensional stochastic approximation process and identify the driving vector field; this will be useful for some of our proofs and specific examples.  Note that zeroes of $P$ correspond to fixed points of $R$, and vice versa; we will describe a zero of $P$ as (linearly) stable, (linearly) unstable or a touchpoint according to the behaviour of the corresponding fixed point of $R$.

Let $\hier t{i,j}$, $i=1,\ldots,N$, $j=1,2$, be the proportion of edge-ends of $G_t$ which are in community $i$ and at a vertex of type $j$; $j=1$ will correspond to red and $j=2$ to blue. Let $\mathbf{X}^{(t)}=\bigl(\hier t{1,1},\hier t{1,2},\ldots,\hier t{N,1},\hier t{N,2}\bigr)$, and note that this takes values in the set \[\Delta^{2N-1}=\biggl\{(x_{1,1},x_{1,2},\ldots,x_{N,1},x_{N,2})\in \bigl(\RR_{\geq 0}\bigr)^{2N}: \sum\nolimits_{i=1}^N (x_{i,1}+x_{i,2})=1\biggr\},\] which is a $(2N-1)$-dimensional set. Write $\hier[Y]ti=\hier t{i,1}+\hier t{i,2}$ for the total proportion of edge-ends in community $i$; note that this means 
$\hier[Y]ti\hier[Z]ti=\hier t{i,1}$.  Write $\salj_t$ for the $\sigma$-algebra generated by all the graphs $G_0,G_1,\ldots,G_t$ and the types of all their vertices.

The new edge-ends added at time $t+1$ which are in community $i$ and at a vertex of type $j$ can be split according to whether they are at the new vertex $v_{t+1}$ or at an existing vertex which it connects to.  For those at $v_{t+1}$, define the random variable $\hier[N]{t+1}{i,j}$ to be equal to $m$ if $c(v_{t+1})=i$ and $v_{t+1}$ has type $j$ and to be equal to zero otherwise.  For those at an existing vertex, define $\hier[E]{t+1}{i,j}$ to be the number of vertices (counted with multiplicity) in community $i$ and of type $j$ which $v_{t+1}$ connects to.

The probability, conditional on $\salj_t$ and that $c(v_{t+1})=i$, that a specific edge from $v_{t+1}$ connects to a vertex of type $j$ is $Q_{i,j}(\mathbf{X}^{(t)})$, where for $x\in\Delta^{2N-1}$ we define \[Q_{i,j}(x)=\frac{\sum_{k=1}^N \alpha_{k,i}x_{k,j}}{\sum_{k=1}^N \alpha_{k,i}(x_{k,1}+x_{k,2})}.\] (Note that $Q_{i,2}(x)=1-Q_{i,1}(x)$.)  Thus the probability that a new vertex at location $i$ at time $t+1$ becomes red (type $1$) is $R(Q_{i,j}(\mathbf{X}^{(t)}))$,
and so \[\EE(\hier[N]{t+1}{i,1}\mid\salj_t)=G_{i,1}(\mathbf{X}^{(t)}),\] where for $x\in\Delta^{2N-1}$ we define 
$G_{i,1}(x)=m\mu_iR\left(Q_{i,1}(x)\right)$.

Similarly we define $G_{i,2}(x)=m\mu_i (1-R(1-Q_{i,2}(x))$, so that the expected number of new edge ends of type $2$ at location $i$ at time $t+1$ coming from the new vertex, conditional on $G_t$, $\EE(\hier[N]{t+1}{i,2}\mid\salj_t)=G_{i,2}(\mathbf{X}^{(t)})$. This uses the fact that the probability a new vertex is blue (type $2$) given that it has $k$ blue neighbours is $1-p_{m-k}$.

Then we have 
\[\EE(\hier[E]{t+1}{i,j}|\salj_t)=H_{i,j}(\mathbf{X}^{(t)}),\] where for $x\in\Delta^{2N-1}$ we define \[H_{i,j}(x)=mx_{i,j}\sum_{k=1}^N \mu_k\frac{ \alpha_{i,k}}{\sum_{\ell=1}^{N}\alpha_{\ell,k}\left(x_{\ell,1}+x_{\ell,2}\right)}.\]

We assumed that $G_0$ has $2mn_0$ edges, and so the total number of edge-ends in $G_t$ is $2m(t+n_0)$.  Thus we have \[2m(t+n_0+1)\hier{t+1}{i,j}-2m(t+n_0)\hier{t}{i,j}=\hier[N]{t+1}{i,j}+\hier[E]{t+1}{i,j},\] and so \[2m(t+n_0+1)\EE(\hier{t+1}{i,j}\mid\salj_t)-2m(t+n_0)\hier{t}{i,j}=G_{i,j}(\mathbf{X}^{(t)})+H_{i,j}(\mathbf{X}^{(t)}).\]

Hence, writing $\gamma_t=(2m(t+n_0+1))^{-1}$, we can write \[X^{(t+1)}_{i,j}-\hier t{i,j}=\gamma_t (F_{i,j}(\mathbf{X}_{t})+\xi^{(t+1)}_{i,j})\] 
where, again for $x\in \Delta^{2N-1}$, \[F_{i,j}(x)=G_{i,j}(x)+H_{i,j}(x)-2m x_{i,j},\] and the $\xi^{(t+1)}_{i,j}$ are noise terms with $\EE(\xi^{(t+1)}_{i,j}\mid\salj_t)=0$, and so $(\mathbf{X}_t)_{t\geq 0}$ is a Robbins-Monro stochastic approximation process as defined in section 4.2 of Bena\"{\i}m \cite{benaim}.  Furthermore, because exactly $2m$ edge ends are added at each time point, each component of the noise term $\xi^{(t+1)}_{i,j}$ is uniformly bounded in modulus by $2m$, which means that we can apply Proposition 4.2 of \cite{benaim} to show that condition A1 of Proposition 4.1 of \cite{benaim} is satisfied.  Because $\Delta^{2N-1}$ is bounded, assumption A2 of the latter result is also satisfied, and so we can use it to conclude that an interpolated version of $(\mathbf{X}_t)_{t\geq 0}$ is an asymptotic pseudotrajectory of the vector field on $\Delta^{2N-1}$ driven by $F$.

From results in \cite{geopref9}, we have the following.  \begin{lemma}\label{lem:nu} For each $i \in \{1,\ldots,N\}$ the limit as $t\to\infty$ of the proportion of edge ends in community $i$, $\hier[Y]ti=\hier{t}{i,1}+\hier{t}{i,2}$, almost surely converges to a particular limit $\nu_i$, where the $\nu_i$ are the unique solution satisfying $\nu_i\in[0,1]$ for all $i$ and $\sum_{i=1}^N \nu_i=1$ to the equations \begin{equation}\label{edgeendlim}\nu_i=\frac12 \mu_i+\frac12 \sum_{j=1}^N\mu_j\frac{\alpha_{i,j}\nu_i}{\sum_{k=1}^N \alpha_{k,j}\nu_k}.\end{equation}\end{lemma}

\begin{proof}See \cite{geopref9}, Proposition 3.1.  That the $\nu$ identified by that proposition satisfies \eqref{edgeendlim} follows from the fact that, in the notation defined in \cite{geopref9} immediately before the statement of Proposition 3.1, \eqref{edgeendlim} is precisely the statement that $G_i(\nu)=0$ and that the $\nu$ found in the proof of the proposition is constructed to satisfy this. \end{proof}  The consequence of this result for our situation is that the limit set of $(\mathbf{X}_t)$ must be contained within the set \[\Delta_{\nu}=\{(x_{1,1},x_{1,2},\ldots,x_{N,1},x_{N,2})\in\Delta^{2N-1}:x_{1,1}+x_{1,2}=\nu_1,\ldots, x_{N,1}+x_{N,2}=\nu_N\}.\]  This will be useful in some examples when it comes to identifying possible limits.

In many cases the equations \eqref{edgeendlim} for the $\nu_i$ do not have a closed form solution, but there are some examples where they do; for example when $\mu_i=\frac{1}{N}$ for all $i$ and the matrix $A$ has suitable symmetries it can be seen that $\nu_i=\frac{1}{N}$ for all $i$.

\subsection{Proof of Theorem \ref{stable}, part 2}

We start by analysing the vector field $F$ to prove the second part of Theorem \ref{stable}, that is that synchronised limits which occur with positive probability can only occur at fixed points of $R$.  By the fact that an interpolated version of $(X_t)_{t\geq 0}$ is an asymptotic pseudotrajectory of the vector field driven by $F$, applying Theorem 5.7 of \cite{benaim} for the case where $L(X)$ is a single point shows that any $x\in\Delta^{2N-1}$ which occurs with positive probability as a limit of $X_t$ as $t\to\infty$ is a stationary point of $F$.

By Lemma \ref{lem:nu} we know that any $x=(x_{1,1},x_{1,2},\ldots,x_{N,1},x_{N,2})\in\Delta^{2N-1}$ which is an element of the limit set of $(\mathbf{X}_t)_{t\geq 0}$ satisfies $x_{i,1}+x_{i,2}=\nu_i$ for each $i\in \{1,\ldots,N\}$.  It follows that for such an $x$ we can write \begin{align*}H_{i,j}(x) &= mx_{i,j}\sum_{k=1}^N \mu_k\frac{ \alpha_{i,k}}{\sum_{\ell=1}^{N}\alpha_{\ell,k}\nu_{\ell}} \\ &= m\frac{x_{i,j}}{\nu_i}\sum_{k=1}^N \mu_k\frac{ \alpha_{i,k}\nu_i}{\sum_{\ell=1}^{N}\alpha_{\ell,k}\nu_{\ell}} \\ &= m\frac{x_{i,j}}{\nu_i}(2\nu_i-\mu_i),\end{align*} using \eqref{edgeendlim} in the last line.

It follows that at such an $x$ \[F_{i,1}(x) = m\mu_i R(Q_{i,1}(x))-mx_{i,1}\frac{\mu_i}{\nu_i},\] implying that if $x$ is a stationary point of $F$ then \[R(Q_{i,1}(x))=\frac{x_{i,1}}{\nu_i}.\]  If the process converges to a synchronised limit, $\hier[Z]{t}{i,1}\to z^*$ for all $i$, then $\hier[X]{t}{i,1}\to \nu_i z^{*}$ for all $i$, so if $x$ is a possible synchronised limit it will have $\frac{x_{i,1}}{\nu_i}=z^{*}$.  If $x$ is of this form then $Q_{i,1}(x)$ is easily seen to be equal to $z^{*}$, so we obtain $R(z^{*})=z^{*}$.  It follows that any positive probability synchronised limit will have $z^{*}$ being a fixed point of $R$.

\subsection{Proof of Theorem \ref{different}}

We may assume $\eps$ is sufficiently small that $R(z)$ has no other fixed points within $2\eps$ of $z_i$, and that $R'(z)<1-\kappa$ in this range, for some $\kappa>0$ and all $i$.

Note that multiplying any column of the matrix $A_1$ by a positive constant has the effect of multiplying the corresponding column of $A_\theta$ by the same constant, and this does not change the model. Consequently, without loss of generality, we may assume that $A_1$ has been normalised so that its diagonal entries are all $1$, and hence the diagonal entries of $A_\theta$ satisfy $\alpha_{i,i}\equiv 1$. We write the off-diagonal entries of $A_1$ as $\tilde{\alpha}_{i,j}$ so that the off-diagonal entries of $A_{\theta}$ are $\alpha_{i,j}=\theta \tilde{\alpha}_{i,j}$. For a stochastic vector $\pi\in \RR^N$ define
\[\hier[q]{\pi}{i,j}=\frac{\alpha_{i,j}\pi_i}{\sum_{k\in S}\alpha_{k,j}\pi_k}.\]
Note that if the current distribution of edge-ends between communities is given by $\pi$ then $\hier[q]{\pi}{i,j}$ represents the probability with which a new vertex chooses a neighbour from community $i$ given that the new vertex is in community $j$.

Let $\hier[R]ti$ (respectively, $\hier[B]ti$) be the number of red (respectively, blue) edge-ends in community $i$ at time $t$. We will define a system of urns to bound the process. While our first set of urns will be enough to prove the weak form of the result alluded to in Remark \ref{rem:weak-different}, we need to iterate the argument in order to show that $\hier[Z]ti$ converges to some limit for each $i$.

\subsubsection{Urns}

Each set of urns that we consider depends on intervals $I_i\subseteq[0,1]$ for each $i\in S$, and a real number $\delta\in (0,\min_i\nu_i/2)$. The urns will be described in terms of other parameters that depend on $\delta$, and an additional parameter $c$ that we shall show only depends on $A_1$ and $\mu$ (although in order to define $c$ we first need to define some other parameters).

For each $i,j\in S$, set
\begin{equation}q_{i,j}^+=\max\{\hier[q]{\pi}{i,j}\mid\|\pi-\nu\|_\infty\leq\delta\}\text{ and } q_{i,j}^-=\min\{\hier[q]{\pi}{i,j}\mid\|\pi-\nu\|_\infty\leq\delta\}.\label{q-pm-ij}\end{equation}

\begin{lemma}\label{lem:thetadelta}There is some constant $c$ depending only on $\mu$ and $A_1$ (in particular, independent of $\delta$ and $\theta$) satisfying
$q^+_{i,j}-q^-_{i,j}\leq c\theta\delta$
for all $i,j\in S$, $\theta>0$ and $\delta\in(0,\min_k\mu_k/4)$,
and
$\hier[q]{\nu}{i,j}\leq c\theta$
for all $i\neq j$ and $\theta>0$. 
\end{lemma}
\begin{proof}
Recall that the limiting distribution of edge-ends in communities, $\nu$, depends on the community structure and hence on the parameter $\theta$.
We first need a lower bound on $\nu_k$ which does not depend on $\theta$. Note that by time $t$ community $k$ receives $\mu_kt-o(t)$ new vertices, contributing at least $\mu_ktm-o(t)$ edge-ends, almost surely, and this makes up a $\mu_k/2-o(1)$ proportion of all edge-ends. Thus $\min_k\nu_k\geq\min_k\mu_k/2$. 

It suffices to prove the two statements separately (and then take the smaller of the two values of $c$). To prove the second, observe that $\alpha_{i,j}\nu_i\leq \theta\max_{k,\ell}\tilde{\alpha}_{k,\ell}$ and, since we are assuming $\alpha_{j,j}=1$, we have $\sum_{k\in S}\alpha_{k,j}\nu_k\geq \nu_j\geq\min_k\mu_k/2$.

To prove the first statement, observe that the maximum is attained when $\pi_i=\nu_i+\delta$ and $\pi_k=\nu_k-\delta$ for $k\neq i$, and conversely for the minimum. Thus
\[q^+_{i,j}-q^-_{i,j}=\frac{2\delta\alpha_{i,j}\sum_{k\neq i}\alpha_{k,j}(\nu_i+\nu_k)}{(\sum_{k\neq i}\alpha_{k,j}(\nu_k-\delta)+\alpha_{i,j}(\nu_i+\delta))(\sum_{k\neq i}\alpha_{k,j}(\nu_k+\delta)-\alpha_{i,j}(\nu_i-\delta))}.\]
Since $\delta<\min_k\mu_k/4\leq\min_k\nu_k/2$, the denominator is at least $(\sum_k\alpha_{k,j}\nu_k/2)^2=a+a'\theta+a''\theta^2$ for some constants $a,a',a''>0$ (which depend on $i$ and $j$). Also, the numerator is of the form $\theta\delta(b+b'\theta)$ for some constants $b>0$ and $b'\geq 0$ (which also depend on $i$ and $j$; $b'>0$ if and only if $i\neq j$). Thus $q^+_{i,j}-q^-_{i,j}<\theta\delta\max\{b/a,b'/a'\}$, and taking $c$ to be the greatest value of $\max\{b/a,b'/a'\}$ as $i,j$ ranges over all pairs gives the required result.
\end{proof}
For each $i$, set \[\xi_i=4c\theta\delta\mu_i^{-1};\]
note that our condition on $\delta$ ensures $\xi_i<c\theta$.

We define urns $\hat U_i$ and $\check U_i$ for each $i\in S$. Denote the number of red balls in urn $\hat U_i$ at time $t$ by $\hier[\hat{R}]ti$, and define $\hier[\hat{B}]ti,\hier[\check{R}]ti,\hier[\check{B}]ti$ similarly. Denote the proportion of red balls in each urn by $\hier[\hat{r}]ti$ or $\hier[\check{r}]ti$ as appropriate.

We initialise the urns at some time $t_{\mathrm{init}}$, with $\hier[\hat{R}]{t_{\mathrm{init}}}i=\hier[\check{R}]{t_{\mathrm{init}}}i=\hier[R]{t_{\mathrm{init}}}i$ and $\hier[\hat{B}]{t_{\mathrm{init}}}i=\hier[\check{B}]{t_{\mathrm{init}}}i=\hier[B]{t_{\mathrm{init}}}i$. The urns are then updated as follows, for each time $t> t_{\mathrm{init}}$. The random variables used are independent of those from previous time steps (but may depend on each other, since coupling each to the graph process will create indirect dependence). 
\begin{enumerate}
\item Let $J(t)$ be a random variable with support $S$ and law $\mu$.
\item For each $i\in S$, do the following.
\begin{itemize}
\item Let $P^{(t)}_i\leq Q^{(t)}_i$ be coupled random variables with $P^{(t)}_i\sim\bin{m,q_{i,J(t)}^-}$ and $Q^{(t)}_i\sim\bin{m,q_{i,J(t)}^+}$.
\item Add $P^{(t)}_i$ balls, which are independently red with probability $\hier[\hat{r}]{t-1}i$ and blue otherwise, to $\hat U_i$, and add an identical group of balls to $\check U_i$.
\item Add $Q^{(t)}_i-P^{(t)}_i$ red balls to $\hat U_i$ and $Q^{(t)}_i-P^{(t)}_i$ blue balls to $\check U_i$.
\end{itemize}
\item Let $r_+^{(t-1)}$ and $r_-^{(t-1)}$ be the maximum value and minimum value respectively of 
$R(\sum_{i\in S}q_ix_i)$, where $x_{J(t)}\in[\hier[\hat{r}]{t-1}{{J(t)}}-\xi_{{J(t)}},\hier[\hat{r}]{t-1}{{J(t)}}+\xi_{{J(t)}}]$ but $x_i\in I_i$ for $i\neq J(t)$, and $q_i\in[q_{i,{J(t)}}^-,q_{i,{J(t)}}^+]$ for each $i$.
\item Let $A_-^{(t)}\leq A_+^{(t)}$ be coupled Bernoulli random variables with parameters $r_-^{(t-1)}$ and $r_+^{(t-1)}$ respectively. Add $m$ balls to $\hat U_{J(t)}$, which are red if $A_+^{(t)}=1$ and otherwise blue, and add $m$ balls to $\check U_{J(t)}$, which are red if $A_-^{(t)}=1$ and otherwise blue.
\end{enumerate}
We next prove that this system of urns may be appropriately coupled to the graph process. 

\subsubsection{Coupling}
We first need the straightforward fact that the number of edge-ends of each type meeting each community tends to infinity.
\begin{lemma}\label{many-blues}Suppose that the initial state has at least one vertex of each colour in each community. Then the total degree of vertices of a given colour in a given community tends to infinity almost surely.
\end{lemma}
\begin{proof}Note that the number of edge-ends meeting a given community $i$ at time $t$ is $\Theta(t)$, and, since $\alpha_{i,i}\mu_i>0$, the probability of any given vertex in community $i$ being selected as a neighbour of $v_t$ is $\Omega(1/t)$, independently of previous selections. Thus, by the second Borel--Cantelli lemma, the degree of this vertex almost surely tends to infinity.
\end{proof}

\begin{lemma}\label{lem:coupling}For any $\zeta>0$ there exists $t_\zeta$ such that if $t_{\mathrm{init}}\geq t_\zeta$ then with probability at least $1-\zeta$ either
\begin{enumerate}[(i)]
\item\label{intervals} $\hier[Z]ti\not\in I_i$ for some $i$ and some $t\geq t_{\mathrm{init}}$, or
\item we may couple the processes with $\hier[\check{r}]ti\leq\hier[Z]ti\leq \hier[\hat{r}]ti$, for all $i\in S$ and $t\geq t_{\mathrm{init}}$.
\end{enumerate}
\end{lemma}
\begin{proof}We assume throughout that
\begin{equation}\abs{\hier[Y]ti-\nu_i}<\delta,\text{ for all }i\text{ and for all }t\geq t_{\mathrm{init}},
	\label{close-to-nu}
\end{equation}
which holds for $t_{\mathrm{init}}$ sufficiently large with probability at least $1-\zeta/3$. We will show by induction the stronger statement that we may couple the two processes such that $\hier[\hat{R}]ti+\hier[\hat{B}]ti=\hier[\check{R}]ti+\hier[\check{B}]ti\geq \hier[R]ti+\hier[B]ti$ for each $i$, and $\hier[\hat{B}]ti\leq \hier[B]ti$ but $\hier[\check{R}]ti\leq \hier[R]ti$. By symmetry, and since it is immediate that $\hat U_i$ and $\check U_i$ have the same number of balls, it suffices to prove the statements for $\hat U_i$.

For $t>t_{\mathrm{init}}$, we consider the steps in the definition of the urn process in turn.
\begin{enumerate}
\item Since they have the same law, we couple $J(t)$ to be the community of $v_t$.
\item For each $i$, let $S^{(t)}_i$ be the number of edges formed from $v_t$ to community $i$. We have $S^{(t)}_i\sim\bin{m,q_i}$, where $q_i=\hier[q]{\pi}{i,J(t)}$ and $\pi=(\hier[Y]{t-1}1,\ldots,\hier[Y]{t-1}N)$. By \eqref{q-pm-ij} and \eqref{close-to-nu}, we have $q_{i,J(t)}^-\leq q_i\leq q_{i,J(t)}^+$ and hence can couple the variables such that $P^{(t)}_i\leq S^{(t)}_i\leq Q^{(t)}_i$. Consequently the number of balls added in this step is at least the number of edge-ends added to existing vertices in community $i$.

Each of the $S^{(t)}_i$ neighbours of $v_t$ in community $i$ is independently red with probability $\hier[Z]{t-1}i\leq\hier[\hat{r}]{t-1}i$ by the induction hypothesis.  Thus we may couple the first $P^{(t)}_i$ balls with the first $P^{(t)}_i$ neighbours such that at most as many of the former are blue. Since all the remaining $Q^{(t)}_i-P^{(t)}_i\geq S^{(t)}_i-P^{(t)}_i$ balls are red, the number of blue balls added is at most the number of edge-ends added to existing blue vertices.
\item Let $p=\prob{v_t\text{ red}\mid c(v_t)=J(t)}$. Then provided $\Bigl|\hier[\hat{r}]{t-1}{J(t)}-\hier[Z]{t-1}{J(t)}\Bigr|\leq\xi_{J(t)}$ and $\hier[Z]{t-1}i\in I_i$ for all $i\neq J(t)$, we have $r_-^{(t)}\leq p\leq r_+^{(t)}$.
\item These $m$ balls correspond to the $m$ new edge-ends meeting $v_t$. If $\hier[\hat{r}]{t-1}{J(t)}\leq \hier[Z]{t-1}{J(t)}+\xi_{J(t)}$, then either \eqref{intervals} holds or $r_+^{(t)}\geq p$, meaning that we may couple the balls to be red whenever $v_t$ is. The coupling of new balls and new edge-ends, together with the induction hypothesis, gives the required inequalities for $t$. 

Furthermore, if $\hier[B]{t-1}{J(t)}-\hier[\hat{B}]{t-1}{J(t)}\geq m$ then, whatever colour these $m$ balls are, we still have $\hier[B]{t}{J(t)}\geq \hier[\hat{B}]{t}{J(t)}$, giving the required inequalities for $t$.
\end{enumerate}
Thus, assuming \eqref{intervals} does not occur, coupling can only fail if $\hier[B]{t-1}{J(t)}-\hier[\hat{B}]{t-1}{J(t)}<m$ but $\hier[\hat{r}]{t-1}{J(t)}-\hier[Z]{t-1}{J(t)}>\xi_{J(t)}$. In this case we have
\begin{align}
	\notag\frac{\hier[\hat{R}]t{J(t)}+\hier[\hat{B}]t{J(t)}}{\hier[R]t{J(t)}+\hier[B]t{J(t)}}\cdot
	\frac{\hier[\hat{B}]t{J(t)}+m}{\hier[\hat{B}]t{J(t)}}&\geq\frac{\hier[\hat{R}]ti+\hier[\hat{B}]ti}{\hier[R]ti+\hier[B]ti}
	\cdot\frac{\hier[B]ti}{\hier[\hat{B}]ti}\\
	&\notag=\frac{1-\hier[Z]{t-1}{J(t)}}{1-\hier[\hat{r}]{t-1}{J(t)}}\\
	&>1+\xi_{J(t)}.\label{too-many-extra}
\end{align}
However, at each step of the process the number of balls added to $\hat U_i$ exceeds the number of edge-ends added to community $i$ by at most $Q^{(t)}_i-P^{(t)}_i\sim\bin{m,q_{i,j}^+-q_{i,j}^-}$ for some $j$, which is at most $\bin{m,c\theta\delta}$, whereas with probability $\mu_i$ we add at least $m$ edge-ends. Thus if $t_{\mathrm{init}}$ is sufficiently large, we have 
\begin{equation}\frac{\hier[\hat{R}]ti+\hier[\hat{B}]ti}{\hier[R]ti+\hier[B]ti}\leq 1+\xi_i/2\text{ for all }i\in S\text{ and }t\geq t_{\mathrm{init}}
	\label{few-extra-balls}\end{equation}
with probability at least $1-\zeta/3$. However, \eqref{few-extra-balls} contradicts \eqref{too-many-extra} provided $\hier[\hat{B}]{t_{\mathrm{init}}}i$ is sufficiently large for each $i$, which also holds with probability at least $1-\zeta/3$ for $t_{\mathrm{init}}$ sufficiently large by Lemma \ref{many-blues}. Thus the failure probability is at most $\zeta$, as required.
\end{proof}

\subsubsection{Bounding}
To get the desired bounds we will need to show that, for $\theta$ and $\delta$ sufficiently small, each urn behaves like a one-community process governed by a polynomial close to $R$. We first give some properties these urns must satisfy.

\begin{lemma}\label{lem:upcrossing}
	Consider an urn process, initialised with $X_0$ balls of which $Y_0$ are red and the rest blue. Let $X_t\leq k$ balls, of which $Y_t$ are red and the rest blue, be added at each time $t\geq 1$, where the $X_t$ and $Y_t$ are random variables. Let $Z_t$ be the proportion of red balls in the urn at time $t$. Suppose that $\mathbb{E}(Y_{t+1}\mid\mathcal{H}_t)\leq f(Z_t)\mathbb{E}(X_{t+1}\mid\mathcal{H}_t)$, where $f$ is a continuous function and $\mathcal{H}_t$ is the filtration $\sigma(X_0,Y_0,\ldots, X_t,Y_t)$. Let $r_1<r_2$ satisfy $f(r)<r$ for all $r\in[r_1,r_2]$. 
	\begin{enumerate}[(i)]
		\item We have $\prob{Z_t\in[r_1,r_2]\text{ for all }t\geq 0}=0$.
		\item Suppose that $Z_0\leq r_1$ and the number of balls at time $0$ is $n$. Then with high probability (as $n\to\infty$) the proportion of red balls never exceeds $r_2$.
	\end{enumerate}
\end{lemma}
\begin{proof}
We prove (ii) first. Note that by skipping any step where $X_t=0$ we may assume that $X_t\geq 1$ and hence the urn has at least $n+t$ balls at each time $t$. 

Fix a time $t_0$ and suppose that $Z_{t_0}\in[r_1,(r_1+r_2)/2]$. Consider a process $Z'_t$ defined by $Z'_{t_0}=Z_{t_0}$ and
\[Z'_{t+1}=\begin{cases} Z_{t+1}\text{ if }Z'_t\in[r_1,r_2]\\
	Z'_t\text{ otherwise}.\end{cases}\]
The conditions ensure that $Z'_t$ is a supermartingale, with increments $O(1/(n+t))$ (where the implicit constant depends only on $k$). Thus, by the Azuma--Hoeffding inequality, $\prob{Z'_{t}>r_2}\leq \exp(-c(n+t_0))$ for some constant $c$ depending only on $k$ and $r_2-r_1$, for any $t\geq t_0$. In particular, since the events $(Z'_t>r_2)$ are nested, \[\prob{\exists t_1\geq t_0:Z_{t_1}>r_2}=\prob{\exists t_1\geq t_0:Z'_{t_1}>r_2}\leq \exp(-c(n+t_0)).\] 

Suppose the proportion of red balls eventually exceeds $r_2$. Then there is some first time $t_1$ when it exceeds $r_2$, and there is some last time $t_0-1$ when it previously was at most $r_1$. If $n$ is sufficiently large it follows that the proportion of red at time $t_0$ is in $[r_1,(r_1+r_2)/2]$. But, by a union bound over possible values of $t_0$, with high probability such a pair $t_0,t_1$ does not exist.

Next we prove (i). By continuity of $f$ there is some $\eps>0$ such that $f(r)\leq r-\eps$ on the interval. Suppose $Z_0\in[r_1,r_2]$ and let $T$ be the (possibly infinite) stopping time $\min\{t:Z_t\not\in[r_1,r_2]\}$. Consider instead the process $Z''_t$ defined by $Z''_{0}=Z_{0}$ and
\[Z''_{t+1}=\begin{cases} Z_{t+1}\text{ if }t<T\\
	Z''_t-\eps/(kt+n)\text{ otherwise}.\end{cases}\]
Now we have $W_t:=Z''_t+\sum_{s<t}\eps/(ks+n)$ is a supermartingale, and Azuma--Hoeffding implies $\prob{W_t<r_2+1}\to 1$ as $t\to\infty$. But for sufficiently large $t$ this event implies $t\geq T$, so $T$ is almost surely finite.
\end{proof}

For each $i\in S$ and $z\in[0,1]$, let $Y_i(z)$ be the set of possible values of $\sum_{j\in S}q_jx_j$ for $x_i\in[z-\xi_i,z+\xi_i]$, $x_j\in I_j$ if $j\neq i$ and $q_j\in[q_{i,j}^-,q_{i,j}^+]$. Let $r^+_i(z)$ and $r^-_i(z)$ be the values of $y\in Y_i(z)$ at which $R(y)$ is maximised and minimised respectively. Now set
\begin{gather*}f^+_i(z)=\frac{z\sum_j\mu_jq_{i,j}^++(1-z)\sum_j\mu_j(q_{i,j}^+-q_{i,j}^-)+R(r^+_i(z))\mu_i}
{\sum_j\mu_jq_{i,j}^++\mu_i}\\
f^-_i(z)=\frac{z\sum_j\mu_jq_{i,j}^+-z\sum_j\mu_j(q_{i,j}^+-q_{i,j}^-)+R(r^-_i(z))\mu_i}
{\sum_j\mu_jq_{i,j}^++\mu_i}.
\end{gather*}
Note that continuity of $R$ implies that both $f^+_i$ and $f^-_i$ are continuous. The ratio of the expected number of red balls added to $\hat{U}_i$ at time $t$ to the expected number of balls added is $f^+_i(\hat{r}_{t-1})$, and $f^-_i$ has the same relationship with $\check{U}_i$.

We thus wish to show that $f^+_i(z)<z$ for $z$ in some suitable interval. Note that $f^+_i(z)$ is a weighted average of $z$ and $g^+_i(z):=R(r^+_i(z))+\mu_i^{-1}(1-z)\sum_j\mu_j(q_{i,j}^+-q_{i,j}^-)$, so it suffices to show that $g^+_i(z)<z$. Similarly we may define $g^-_i(z)$.

\begin{lemma}\label{urn-bounds}There is some constant $c'$, which depends only on $A_1$, $\mu$ and $R$, such that the following holds. Suppose that $\abs{I_j}\leq h$ for each $j$. Then for any $z,z'\in[0,1]$ and any $i\in S$ we have $\abs{g^+_i(z)-R(z)},\abs{g^-_i(z)-R(z)}\leq c'\theta$ and $\abs{r^+_i(z)-r^-_i(z')}\leq\abs{z-z'}+c'\theta(h+\delta)$.
\end{lemma}
\begin{proof}Suppose $y=\sum_{j\in S}q_jx_j\in Y_i(z)$ and $y'=\sum_{j\in S}q_jx'_j\in Y_i(z')$. Then
\begin{align*}\abs{y-y'}&\leq\sum_j(x_j\abs{q_j-q'_j}+q'_j\abs{x_j-x'_j})\\
	&\leq Nc\theta\delta+\abs{z-z'}+2\xi_i+(N-1)c\theta(1+\delta) h,
\end{align*}
using the bounds in Lemma \ref{lem:thetadelta}, and that $q'_j\leq q^{(\nu)}_{i,j}+\abs{q^+_{i,j}-q^-_{i,j}}\leq c\theta(1+\delta)$ for $j\neq i$.
By definition of $\xi_i$, this is $\abs{z-z'}+O(\theta\delta+\theta h)$.	
	
Note that $z=\sum_{j\in S}\hier[q]{\nu}{i,j}z$ and hence, using Lemma \ref{lem:thetadelta}, 
\begin{align*}\abs{y-z}&\leq\sum_j(x_j\abs{q_j-\hier[q]{\nu}{i,j}}+\hier[q]{\nu}{i,j}\abs{x_j-z})\\
&\leq Nc\theta\delta+\xi_i+(N-1)c\theta.
\end{align*}
Recalling that $\xi_i\leq c\theta$, and $\delta\leq\min_j\mu_j/4$, it follows that $r^+_i(z)-z=O(\theta)$. Since $R$ is a polynomial, its derivative is uniformly bounded on $[0,1]$, and thus we also have $\abs{R(r^+_i(z))-R(z)}=O(\theta)$ (and similarly for $r^-_i(z)$). The required bounds follow since $0\leq \sum_j\mu_j(q_{i,j}^+-q_{i,j}^-)\leq c\theta\delta$ by Lemma \ref{lem:thetadelta}.
\end{proof}

Now we describe how to set up the urns. We first set $\delta=\delta_0$ and $\zeta=\zeta_0$ sufficiently small (to be chosen later). Set up a system of urns with $I_i=[0,1]$ for each $i$. By Lemma \ref{lem:coupling}, since \eqref{intervals} cannot occur, we may initialise these urns at any sufficiently large time $t_0$ and couple them thereafter with probability at least $1-\zeta_0$.

Recall that $R(z_i+\eps)<z_i+\eps$. It follows from Lemma \ref{urn-bounds} and continuity that for $\theta$ sufficiently small we have $g^+_i(z)<z$ for $z$ in some neighbourhood of $z_i+\eps$. Consequently, if $t_0$ is sufficiently large and $\hat{r}_{t_0}$ is sufficiently close to $z_i$ then Lemma \ref{lem:upcrossing} (ii) ensures that with probability at least $1-\zeta_0/N$ we have $\hier[\hat{r}]ti<z_i+\eps$ for all $t\geq t_0$. Similarly, under the same initial conditions we have $\hier[\check{r}]ti>z_i-\eps$ for all $t\geq t_0$ with probability at least $1-\zeta_0/N$. By applying the same argument to each community, with probability at least $1-3\zeta_0$ each remains within $\eps$ of its desired fixed point for all sufficiently large time.

Thus far we have not used linear stability, and this proves the weak form mentioned in Remark \ref{rem:weak-different}. If we knew that each $\hier[Z]ti$ converged to some limit, we would be done. The remainder of the proof shows that this occurs with positive probability given the bounds above.

To complete the proof, we now iteratively construct a sequence of urn processes. In iteration $n$, first write $\hier[\hat{s}]ti,\hier[\check{s}]ti$ for the proportions of red in the urns from iteration $n-1$. We choose $\delta_n>0$ and $\zeta_n>0$ and take $I_i$ to be the interval in which $\hier[\hat{s}]ti,\hier[\check{s}]ti$ is already known to lie for $t$ sufficiently large, \ie for iteration $1$ we take $I_i=[z_i-\eps,z_i+\eps]$. We ensure that the start time $t_n$ is sufficiently large that $\hier[\hat{s}]ti,\hier[\check{s}]ti\in I_i$ for all $i$ and all $t\geq t_n$ with probability $1-\zeta_n$.

By Lemma \ref{lem:coupling}, and the choice of $I_i$, assuming the success of all previous iterations we may couple this set of urns to the communities successfully with probability at least $1-3\zeta_n$ for a sufficiently large start time $t_n$. Note also that there is a natural coupling of each urn to the corresponding urn from the previous round, which ensures that $\hier[\hat{r}]ti\leq\hier[\hat{s}]ti$ and $\hier[\check{r}]ti\geq\hier[\check{s}]ti$. Providing both couplings to $\hier[\hat r]ti$ succeed, giving $\hier[\hat{r}]ti\geq\hier[\check{r}]ti$, this ensures $\hier[\hat{r}]ti,\hier[\check{r}]ti\in I_i$.

We next define an interval $J_i\subset I_i$ such that the above implies that $\hier[\hat{r}]ti,\hier[\check{r}]ti\in J_i$ for all $t$ sufficiently large almost surely. Indeed, write $I_i=[z_L,z_U]$ and define $z'_L=\inf\{z\in I_i:g^-_i(z)\leq z\}$ and $z'_U=\sup\{z\in I_i:g^+_i\geq z\}$. Then $J_i=[z'_L,z'_U]$ has the required properties, since Lemma \ref{lem:upcrossing} ensures that $\hier[\hat{r}]ti$ crosses $[z'_U,(z_U+z'_U)/2]$ (in that direction) only finitely many times, but also that after each crossing it eventually leaves $[z'_U,z_U]$ (the lower bound on $\hier[\check{r}]ti$ is similar).

Take $h$ minimal such that each $I_i$ has width at most $h$, and choose $\delta_n\leq h$. By continuity we have $g^+(z'_U)\geq z'_U$ and $g^-_i(z'_L)\leq z'_L$. Thus 
\[R(r^+_i(z))-R(r^-_i(z))\geq z'_U-z'_L-2\mu_i^{-1}c\theta h.\]
Also, by Lemma \ref{urn-bounds} we have $r^+_i(z)-r^-_i(z)\leq z'_U-z'_L+2c'\theta h$.
Thus, we obtain by the mean-value theorem that
\begin{align*}R'(z)&\geq\frac{z'_U-z'_L-2\mu_i^{-1}c\theta h}{z'_U-z'_L+2c'\theta h}\\
	&>1-\frac{2(\mu_i^{-1}c+c')\theta h}{z'_U-z'_L}\end{align*}
for some $z\in(z_L,z_U)$.

Since we have $R'(z)<1-\kappa$ in the required range, for some $\kappa>0$ which does not depend on $\theta$, we obtain $z'_U-z'_L<2(\mu_i^{-1}c+c')\theta h/\kappa$. For $\theta$ sufficiently small (independent of $h$), we have $\abs{J_i}\leq \lambda h$ for all $i$, where $\lambda<1$ is constant.

Iterating this procedure therefore gives upper and lower bounds for community $i$ which converge to the same value, and this holds for every community with probability at least $1-3\sum_{n\geq 0}\zeta_n>0$ by suitable choice of $\zeta_n$.

\subsection{Proofs of Proposition \ref{touchpoint}, Theorem \ref{same1}, Theorem \ref{same2}, Theorem \ref{determinant} and Theorem \ref{stable}, part 1}
We shall use the following result more than once.
\begin{lemma}\label{lem:shift-up}Let $z_1<z_2$ be such that $R(z)>z_1$ for $z\in[z_1,z_2)$, and either $R(z_2)>z_1$ or $\alpha_{ii}>0$ for each $i$. Then there exists $\eps>0$ such that almost surely one of the following holds:
\begin{enumerate}[(i)]\item $\lim\inf\hier[Z]ti<z_1-\eps$ for some $i$;
\item $\lim\sup\hier[Z]ti>z_2+\eps$ for some $i$; or
\item $\lim\inf\hier[Z]ti>z_1+\eps$ for every $i$.
\end{enumerate}
\end{lemma}
\begin{figure}[!htb]
	\begin{tikzpicture}
		\begin{axis}[
			width=.9\textwidth,
			height=8cm,
			axis lines = left,
			ymin=0,ymax=0.6,
			xlabel = $z$,
			xmin=0, xmax=1,
			ylabel = $R(z)$,
			]
			\addplot [
			domain=0:1, 
			samples=300, 
			color=black,
			thick,
			]
			{-4.5*x^5 + 20*x^4 - 25*x^3 + 10*x^2};
			\addplot [
			domain=0:0.6, 
			samples=300, 
			color=black,
			]
			{x};
			\addplot [
			thick,
			color=red,
			fill=red,
			fill opacity=0.05]
			coordinates {
			(0.4, 0) 
			(0.4, 0.4)
			(0.6, 0.4)
			(0.6,0)  };
			\addplot [
			thick,
			color=black,
			dotted,
			]
			coordinates {
			(0.39, 0) 
			(0.39, 0.425)
			(0.61, 0.425)
			(0.61,0)  };
		\end{axis}
	\end{tikzpicture}
	\caption{An example of $R(z)$ satisfying the conditions of Lemma \ref{lem:shift-up} with $z_1=0.4$ and $z_2=0.6$. The curve must lie outside the highlighted region. For some suitable $\eps>0$ the value of $R(z)$ is above $z_1+4\mu_i^{-1}\eps$ for $z\in(z_1-2\eps,z_2+2\eps)$, as shown by the dotted lines.}\label{fig:shift-up}
\end{figure}

\begin{proof}Suppose that with positive probability (i) and (ii) do not occur, for some value of $\eps$ to be chosen later. We condition on this event and show that (iii) almost surely occurs. It is sufficient to show that for an arbitrary $i$ and suitable $\eps=\eps_i$ we obtain $\lim\inf\hier[Z]ti>z_1+\eps$, since then (iii) follows for $\eps=\min_i\eps_i$.

Recall that at time $t$ community $i$ gains at most $m$ edge-ends meeting existing vertices, and gains $m$ edge-ends meeting a new vertex with probability $\mu_i$. The former are each independently red with probability $\hier[Z]ti$, and the latter are red with probability $R(r_t)$, where 
$r_t=\frac{\sum_j\alpha_{ji}\hier[Y]tj\hier[Z]tj}{\sum_j\alpha_{ji}\hier[Y]tj}$.
In particular, since $r_t$ is a convex combination of $\hier[Z]tj$ we have $\min_j\hier[Z]tj\leq r_t\leq\max_j\hier[Z]tj$. Thus the expected number of edge-ends added to community $i$ at time $t$ is $e_t=m\mu_i+k_t$, for some $k_t\leq m$, and the expected number of red edge-ends amongst these is $f_t=m\mu_iR(r_t)+k_t\hier[Z]ti$. We will show that almost surely for all sufficiently large $t$ we have $f_t> (z_1+\eps)e_t$; since the number of edges added at each time is bounded, (iii) follows. 

Note that $f_t/e_t$ is monotonic in $k_t$, being increasing if $\hier[Z]ti>R(r_t)$ and decreasing otherwise. Thus, since $0\leq k_t\leq m$, we have \begin{equation}\label{f-e-bounds}f_t/e_t\geq\min\left\{ R(r_t),\frac{\mu_iR(r_t)+\hier[Z]ti}{\mu_i+1}\right\}.\end{equation}
Suppose $R(z_2)>z_1$. Then by continuity of $R$, there is some $\eps>0$ such that $R(z)>z_1+4\mu_i^{-1}\eps$ for $z\in[z_1-2\eps,z_2+2\eps]$ (see Figure \ref{fig:shift-up}. For all sufficiently large $t$ we have $\min_j\hier[Z]tj>z_1-2\eps$ and $\max_j\hier[Z]tj<z_2+2\eps$, so $R(r_t)\geq z_1+4\mu_i^{-1}\eps$. Thus the bound in \eqref{f-e-bounds} is at least $\frac{\mu_i(z_1+4\mu_i^{-1}\eps)+z_1-2\eps}{\mu_i+1}=z_1+\frac{2\eps}{\mu_i+1}>z_1+\eps$, as required.

Now we deal with the case $R(z_2)=z_1$, where we have $\alpha_{ii}>0$. Recall that $\hier[Y]ti\to\nu_i$ almost surely, and so for sufficiently large $t$ we have $\hier[Y]ti>\nu_i/2$. Set $x=\frac{\alpha_{ii}\nu_i}{2\max_j\alpha_{ji}}$. It follows that, for $t$ sufficiently large, we have 
\begin{equation}\label{r-bound}r_t<x\hier[Z]ti+(1-x)\max_j\hier[Z]tj\end{equation}
and $k_t\geq\mu_ixm$. The latter implies
\begin{equation}\label{f-e-bounds2}f_t/e_t\geq\min\left\{\frac{R(r_t)+x\hier[Z]ti}{1+x},\frac{\mu_iR(r_t)+\hier[Z]ti}{\mu_i+1}\right\}.\end{equation}
Let $z_3=z_1x/2+z_2(1-x/2)$; since $z_1<z_3<z_2$ we may find $\eta>0$ such that $R(z)>4\mu_i^{-1}\eta$ for $z\in[z_1-2\eta,z_3+2\eta]$. Choose $\zeta>0$ such that $\frac{x(z_2-z_1)/2-\zeta}{x+1}>\zeta$, and choose $\xi>0$ such that $R(z)>z_1-\zeta$ for $z\in[z_2,z_2+2\xi]$. Finally, set $\eps=\min\{\eta,\zeta,\xi\}$.

For all sufficiently large $t$ we have $\min_j\hier[Z]tj>z_1-2\eps$ and $\max_j\hier[Z]tj>z_2+2\eps$. Provided $\hier[Z]ti<(z_1+z_2)/2$, \eqref{r-bound} gives $r_t<z_3+2\eta$ and hence $R(r_t)\geq z_1+4\mu_i^{-1}\eta$, and as before it follows that $f_t/e_t>z_1+\eta$. If $\hier[Z]ti<(z_1+z_2)/2$ then \eqref{f-e-bounds2} gives $f_t/e_t> \frac{z_1-\zeta+x(z_1+z_2)/2}{1+x}>z_1+\zeta$. Since $\eps\leq\eta,\zeta$ this completes the proof.
\end{proof}

A technical difficulty in some cases is to rule out the possibility of unstable fixed points at $0$ and/or $1$ becoming accumulation points. We give a general result on non-convergence to $0$.
\begin{lemma}\label{lem:nozero}Suppose that either of the following conditions holds:
\begin{enumerate}[(a)]
\item\label{nozero} $R(0)>0$, and either $R(1)>0$ or $\alpha_{ii}>0$ for each $i$;
\item\label{unstable} $R$ has a linearly unstable zero at $0$.
\end{enumerate}
Then there exists some $\eta>0$ such that for each $i$ we have $\lim\inf_{t\to\infty}\hier[Z]ti\geq\eta$ almost surely.
\end{lemma}
\begin{remark}The additional condition on the matrix in \eqref{nozero} is necessary to exclude situations such as in Theorem \ref{same2}, where it is possible for some communities to converge to $0$ and others to $1$ despite neither being a fixed point of $R$.
\end{remark}
\begin{proof}
Note that, since $R(z)\not\equiv 0$, we have $R(z)>0$ for all $z\in (0,1)$. If \eqref{nozero} holds then taking $(z_1,z_2)=(0,1)$ in Lemma \ref{lem:shift-up} immediately gives the desired conclusion.

Now assume \eqref{unstable} holds. If $R(1)=0$ then, by \eqref{nozero} applied after swapping colours, almost surely there exists some $\eta'>0$ such that $\hier[Z]ti\leq 1-\eta'$ for all $i$ and all $t$ sufficiently large. If $R(1)>0$ then we may instead take $\eta'=0$. In either case, since $0$ is a linearly unstable root and $R(1-\eta')>0$, there is some $\lambda>0$ and $\delta_0>0$ such that $R(z)>\min\{(1+\lambda)z,\delta_0\}$ whenever $0<z\leq 1-\eta'$. 

Write $c=\min_{i\in [N]}\mu_i/2$ and $\kappa=\min_i\nu_i/2$. Set $\lambda'=\lambda/(2+c+\lambda c)<1/c$. Choose $\eta=\delta_0\frac{1-c\lambda'}{1+2\lambda'}$.

For each $t$, set $W_t=\min_{i\in[N]}\hier[Z]ti$. We choose a sequence of values $\delta_1,\delta_2,\ldots$, as follows. Choose $t_0$ sufficiently large that, for all $t\geq t_0$, we have $\hier[Z]ti\leq 1-\eta'$ and $\hier[Y]ti\geq\kappa$ for all $i$. Let $t_1\geq t_0$ be the next time satisfying $W_{t_1}<\frac{\delta_0}{1+2\lambda'}$, if it exists, and set $\delta_1=W_{t_1}$. Now, for each $i\geq 2$, let $t'_{i}$ be the next time after $t_{i-1}$ that either $W_{t'_i}\leq\delta_{i-1}(1-c\lambda')$ or $W_{t'_i}\geq\delta_{i-1}(1+c\lambda'/2)$ (if it exists; we shall subsequently argue that one of these events occurs almost surely). Choose $t_i\geq t'_i$ to be the next time with $W_{t_i}<\frac{\delta_0}{1+2\lambda'}$, if it exists, and set $\delta_{i}=W_{t_i}$. 

Suppose $t_{i-1}$ exists for some $i\geq 2$. Note that we have $(1+2\lambda')\delta_{i-1}<\delta_0$, and hence $R(z)\geq(1+2\lambda')\delta_{i-1}$ for all $z\in[(1-c\lambda')\delta_{i-1},1-\eta']$. Thus, while we have $W_t\geq\delta_{i-1}(1-c\lambda')$, each additional edge-end in community $j$ created at a new vertex has probability at least $\delta_{i-1}(1+2\lambda')$ of being red, and each additional edge-end created at an existing vertex has probability at least $\delta_{i-1}(1-c\lambda')$ of being red, so the expected proportion of red edge-ends created in community $j$ is at least $\delta_{i-1}(1+c\lambda')$. Consequently, almost surely eventually either $W_t<\delta_{i-1}(1-c\lambda')$ or $W_t>\delta_{i-1}(1+c\lambda')$, \ie $t'_{i}$ exists. If $t'_{i}$ exists but $t_{i}$ does not, then $W_t\geq\frac{\delta_0}{1+2\lambda'}>\eta$ for all $t\geq t'_{i}$. Thus we may assume $t_i$ exists for every $i$.

We say that there is a ``decrease'' at $i$ if $\delta_i\leq\delta_{i-1}(1-c\lambda')$. Suppose there is no decrease at $i$. Then $W_{t'_i}\geq \delta_{i-1}(1+c\lambda'/2)$. If $W_{t'_i}<\frac{\delta_0}{1+2\lambda'}$ then $T_i=t'_i$, and otherwise we have $W_{t_i}<\frac{\delta_0}{1+2\lambda'}$ but $W_{t_i-1}\geq\frac{\delta_0}{1+2\lambda'}$. Since $|W_{t}-W_{t-1}|=O(1/t)$ it follows that  $\delta_i\geq\min\{\delta_{i-1}(1+c\lambda'/2),\frac{\delta_0}{1+2\lambda'}-o(1)\}$ as $i\to\infty$. Consequently, if there are only finitely many decreases then we have $\lim\inf_{i\to\infty}\delta_i\geq\frac{\delta_0}{1+2\lambda'}$, and therefore $\lim\inf_{t\to\infty}W_t\geq\eta$. We will show that in fact this almost surely happens, giving a contradiction.

We bound the probability of a decrease, \ie that any community reaches $\delta_{i-1}(1-c\lambda')$ before all communities are at least $\delta_{i-1}(1+c\lambda'/2)$, using Hoeffding's inequality. We take a union bound over the following events:
\begin{enumerate}[(i)]
\item for some $t$ with $t_{i-1}<t<3t_{i-1}$, and some $j$, community $j$ is below $\delta_{i-1}(1-c\lambda')$ at time $t$, and $t$ is minimal for this property;
\item none of the above events occur, but for some $j$ community $j$ is below $\delta_{i-1}(1+c\lambda'/2)$ at time $3t_{i-1}$. 
\end{enumerate}
Since there are at least $2mt_{i-1}\kappa$ edge-ends in community $j$ at time $t_{i-1}$, at most $2m(t-t_{i-1})$ edge-ends are added between times $t_{i-1}$ and $t$, and at least a $\delta_{i-1}$ proportion of the former are red, in order for (i) to occur we must have
\[\frac{t_1\kappa\delta_{i-1}}{t_{i-1}\kappa+t-t_{i-1}}\leq\delta_{i-1}(1-c\lambda'),\]
\ie $t-t_{i-1}\geq \rho t_{i-1}$, where $\rho=c\lambda'\kappa/(1-c\lambda')>0$ does not depend on $t_{i-1}$ or $\delta_{i-1}$.

Fix a community $j$ and reveal for each time step between $t_{i-1}$ and $3t_{i-1}$ the numbers of edge-ends added to community $j$ at old and new vertices. The probability that, for some $t$ with $(1+\rho) t_{i-1} <t\leq 3t_{i-1}$, we have had fewer than $0.9\mu_j(t-t_{i-1})$ new vertices added to community $j$ between times $t_{i-1}$ and $t$ is exponentially small in $t_{i-1}$. Assuming this does not occur, and assuming no bad event of type (i) occurs in any community before time $t$, each of these new vertices is red with probability at least $\delta_{i-1}(1+2\lambda')$, and the (at most) $m(t-t_{i-1})$ other edge-ends created in community $j$ are red with probability at least $\delta_{i-1}(1-c\lambda')$, and so we can couple these to binomial distributions of the appropriate probability. For an event of either type to occur, the actual proportion of red edge-ends added between times $t_{i-1}$ and $t$ (or $t_{i-1}$ and $3t_{i-1}$) must be at most $\delta_{i-1}(1+c\lambda'/2)$, requiring one of the two binomial variables to differ from its expectation by a constant multiplicative factor. By a standard Chernoff bound, this has probability which decays exponentially in the expectation, \ie exponentially in $\delta_{i-1}t_{i-1}$. 

Note that, since $t_i\delta_i$ is (up to a constant) the number of red edge-ends in some community, since the number of red edge-ends in each community increases at least logarithmically over time almost surely, and since we must have $t_i\geq t_{i-1}(1+\Theta(1))$ in order for the proportion of red edge-ends in some community to change by a factor of $1-c\lambda'$ or $1+c\lambda'/2$, these probabilities are almost surely summable, and hence by Borel--Cantelli almost surely there are only finitely many decreases, as required.
\end{proof}

\subsubsection{Proof of Proposition \ref{touchpoint}}
Note that we have $R(z)\geq z$ for all $z\leq z^*$, and $R(z)\leq z$ for all $z\geq z^*$. Write $\delta=\min_i\nu_i/2$.

Suppose that with positive probability all communities converge to limits, but that not all limits are equal. Let the random variables $Z_1$ and $Z_2$ be the largest and smallest limits, and let $z_1<z_2$ be values chosen so that, for any $\eps>0$, with positive probability we have $\abs{Z_i-z_i}<\eps/2$ for $i=1,2$. We condition on this event, where $\eps$ is to be chosen later. We may assume without loss of generality that $\frac{z_1+z_2}{2}\leq z^*$, and that community $1$ converges to $Z_1$ but community $2$ converges to $Z_2$. For $t$ sufficiently large each community is within $\eps/2$ of its eventual limit, and the relative size of each community is within $\delta$ of its limit $\nu_i$, for any given $\delta>0$. Note that there exists some $\theta^*>0$, which depends only on the community structure, such that for $\theta<\theta^*$ the probability of a new vertex in any community selecting a neighbour from its own community is at least $2/3$. Furthermore, since each entry of $A_1$ is positive, for each $\theta>0$ there is some $\eta=\eta(\theta)>0$ such that the probability of a new vertex in any community $i$ selecting a neighbour in community $j$ is at least $\eta$ for each pair $i\neq j$.

Consider the probability that a new vertex joining community $1$ becomes red. This is $R(r)$, where $r$ is the probability that the first edge formed by the new vertex is to a red vertex. If $\theta<\theta^*$ and $\eps$ is sufficiently small (in terms of $z_1$, $z_2$ and $\eta$), and $t$ is sufficiently large, then we have $r\leq 2(z_1+\eps)/3+(z_2+\eps)/3\leq z^*$, but also $r\geq (1-\eta)(z_1-\eps)+\eta(z_2-\eps)\geq z_1+2\eps$. Thus $R(r)\geq r\ge z_1+2\eps$, and so almost surely the proportion of red edge-ends in community $1$ will eventually reach $z_1+2\eps$, contradicting our assumptions.

\subsubsection{Proof of Theorem \ref{same1}, and Theorem \ref{stable}, part 1}
First we will express the condition on the fixed points of $R(R(z))$ in Theorem \ref{same1} into an equivalent, but more convenient form.
\begin{lemma}\label{rcond}Suppose that $R(z)$ has a unique stable fixed point $z^*$ in $[0,1]$, and no other fixed points in $(0,1)$.
Then the following are equivalent:
\begin{enumerate}[(i)]\item\label{root-condition} The only fixed points of $R(R(z))$ in $[0,1]$ are the fixed points of $R(z)$;
\item\label{two-point-condition}there do not exist $0\leq z'<z^*<z''\leq 1$ such that $R(z')\geq z''$ and $R(z'')\leq z'$.
\end{enumerate}
\end{lemma}
\begin{figure}[!htb]
	\begin{tikzpicture}
		\begin{axis}[
			width=.9\textwidth,
			height=8cm,
			axis lines = left,
			ymin=0,ymax=1,
			xlabel = $z$,
			xmin=0, xmax=1,
			ylabel = $R(z)$,
			]
			\addplot [
			domain=0:1, 
			samples=300, 
			color=black,
			thick,
			]
			{2.2*x^3 - 3.3*x^2 + 0.3*x + 0.9};
			\addplot [
			domain=0:1, 
			samples=300, 
			color=red,
			thick,
			dashed
			]
			{23.4256*x^9 - 105.415*x^8 + 167.706* x^7 - 95.0334* x^6 - 15.4638* x^5 + 30.2742* x^4 - 5.8146* x^3 + 1.2078* x^2 - 0.0882* x + 0.1008};
			\addplot [
			domain=0:1, 
			samples=300, 
			color=black,
			]
			{x};
			\addplot [
			thick,
			color=black,
			dotted,
			]
			coordinates {
			(0.2, 0) 
			(0.2, 0.8456)
			(0.8, 0.8456)
			(0.8, 0)  };
			\addplot [
			thick,
			color=black,
			dotted,
			]
			coordinates {
			(0.2, 0.1544) 
			(0.8, 0.1544)  };
			\addplot [
			thick,
			color=black,
			dotted,
			]
			coordinates {
			(0.1011, 0) 
			(0.1011, 0.1011)  };							
			\addplot [
			thick,
			color=black,
			dotted,
			]
			coordinates {
			(0.8989, 0) 
			(0.8989, 0.8989)  };
		\end{axis}
	\end{tikzpicture}
	\caption{An example where neither condition in Lemma \ref{rcond} is satisfied; the solid black line shows $R(z)$. The values $z'=0.2$ and $z''=0.8$ violate \eqref{two-point-condition}, and $R(R(z))$ (dashed red line) has linearly stable fixed points close to $0.1$ and $0.9$.}\label{fig:rcond}
\end{figure}
\begin{proof}First we show the forward implication. Suppose \eqref{two-point-condition} is not true; note that this implies $z^*\in(0,1)$. See Figure \ref{fig:rcond} for an example. The set of pairs $(z',z'')$ with $0\leq z'\leq z^*\leq z''\leq 1$ satisfying $R(z')\geq z''$ and $R(z'')\leq z'$ is closed, so there is a pair that maximises $z''-z'$. By assumption, this pair must have $z''-z'>0$, and it follows that $z'<z^*<z''$. We must have either $R(z')=z''$ or $R(z'')=z'$, since otherwise for sufficiently small $\eps>0$ the values $z'-\eps,z''+\eps$ would also violate \eqref{two-point-condition}. Suppose, without loss of generality, that $R(z')=z''$. Then $R(R(z'))=R(z'')\leq z'$. Since $0$ is not a stable fixed point, we have $R(R(\eps))>\eps$ for all sufficiently small $\eps>0$. Consequently there is a fixed point of $R(R(z))$ in $(0,z']$, but this interval does not include $z^*$, and so \eqref{root-condition} is not true.
	
Conversely, suppose that \eqref{root-condition} is not true, and let $z'$ satisfy $R(R(z'))=z'$ but $R(z')\neq z'$. Set $z''=R(z')$; since also $R(R(z''))=z''$ and $z'=R(z'')$, we may assume without loss of generality that $z'< z''$. Thus $R(z')>z'$ and $R(z'')<z''$, and so there is a fixed point of $R$ in the interval $(z',z'')$. Since the only fixed point in $(0,1)$ is $z^*$, we must have $z'<z^*<z''$, and so \eqref{two-point-condition} fails.
\end{proof}

\begin{proof}[Proof of Theorem \ref{same1}]
We have \eqref{two-point-condition} of Lemma \ref{rcond}. We also straightforwardly have $R(z)>z$ for $0<z<z^*$ and $R(z)<z$ for $z^*<z<1$. Write $c=\min_{i\in [N]}\mu_i/2$, and set
\begin{align*}Z_1&=\min_{i\in [N]}\liminf_{t\to\infty}\hier[Z]ti;\\
Z_2&=\max_{i\in [N]}\limsup_{t\to\infty}\hier[Z]ti.
\end{align*}
Suppose for a contradiction that with positive probability these are not both equal to $z^*$, and choose values $(z_1,z_2)\neq(z^*,z^*)$ such that for any $\eps>0$ there is positive probability that $|Z_1-z_1|,|Z_2-z_2|<\eps$. Either these are on the same side of $z^*$, \ie without loss of generality we have $z_1<z^*$ and $z_2\leq z^*$, or they are on opposite sides, \ie $z_1<z^*<z_2$. In the former case, clearly we have $R(z)\geq z>z_1$ for $z_1< z\leq z^*$, and so for all $z\in(z_1,z_2]$. In the latter case, we cannot have both a value $z_3\in[z_1,z^*)$ with $R(z_3)\geq z_2$ and a value $z_4\in(z^*,z_2]$ with $R(z_4)\leq z_1$, since this would imply $R(z_3)\geq z_2\geq z_4$ and $R(z_4)\leq z_1\leq z_3$, contradicting \eqref{two-point-condition} of Lemma \ref{rcond}. Consequently, at least one of these does not exist; assume without loss of generality we have $R(z)>z_1$ for all $z\in(z^*,z_2]$. Then we again have $R(z)>z_1$ for all $z\in(z_1,z_2]$. We also have $R(z_1)>z_1$, since this holds unless both $z_1=0$ and $0$ is a linearly unstable fixed point, which is impossible by Lemma \ref{lem:nozero}. Now Lemma \ref{lem:shift-up} applies, contradicting the assumed properties of $(z_1,z_2)$.
\end{proof}

\begin{proof}[Proof of Theorem \ref{stable}, part 1]Let $(z_1,z_2)$, with $z_1<z^{*}<z_2$, be an interval strictly contained in the neighbourhood of $z^{*}$ for which $R$ is increasing.  Then, setting $c=\min_{i\in[N]}\mu_i/2$, there exists $\eps>0$ such that if $\hier[Z]ti\in (z_1-c\eps,z_2+c\eps)$ for all $i$  the probability of a new edge end in community $i$ created at a new vertex being red is at least $z_1+2\eps$, and the probability of a new edge end in community $i$ at an existing vertex being red is at least $z_1-c\eps$, with an analogous argument showing that these are at most $z_2-2\eps$ and $z_2+c\eps$ respectively.  By choice of $c$ it follows that, for sufficiently large $t_0$, conditional on $\hier[Z]{t_0}i\in (z_1,z_2)$ for all $i$ there will be positive probability that $\hier[Z]ti\in (z_1-c\eps,z_2+c\eps)$ for all $i$ and for all $t\geq t_0$, and on this event $\hier[Z]ti \to z^{*}$ almost surely.  To complete the proof, we observe that for sufficiently large $t_0$ there will always be a choice of the initial steps up to time $t_0$ which has positive probability and gives $\hier[Z]{t_0}i\in (z_1,z_2)$ for all $i$.
\end{proof}

\subsubsection{Proof of Theorem \ref{determinant}, part 1}
We first exclude a degenerate special case. Suppose $\alpha_{2,1}=\alpha_{2,2}=0$, which means that vertices from community $2$ are never selected as neighbours. In this case, community $1$ follows a stochastic approximation equation very similar to that of the single-community case with the same attachment rule, with the only difference being some additional zero-expectation noise when new vertices are added to community $2$. Thus the methods of \cite{AMR} give almost sure convergence to $z^*$ in community $1$. Note that $m$ new edge-ends are added to community $2$ at time $t$ with probability $\mu_2$, and these are red with probability $R(\hier[Z]t1)$. Since $z^*$ is a fixed point of $R$, and $R$ is continuous, almost sure convergence of $\hier[Z]t1$ to $z^*$ implies almost sure convergence of this probability, and hence of $\hier[Z]t2$, to $z^*$.

Thus we may henceforth assume that neither row of $A$ is zero, and since the determinant is non-negative (and no column is zero by definition of the model) it follows that $\alpha_{1,1}$ and $\alpha_{2,2}$ are both positive. Since $0$ either is not a fixed point or is linearly unstable, the conditions of Lemma \ref{lem:nozero} are met.

Define the parameters $Z_i^-=\lim\inf_{t\to\infty}\hier[Z]ti$ and $Z_i^+=\lim\sup_{t\to\infty}\hier[Z]ti$, for $i=1,2$. Suppose for a contradiction that with positive probability these are not all equal to $z^*$, and let $z_1^-,z_1^+,z_2^-,z_2^+$ be specific values such that for any $\eps>0$ there is positive probability that each limit is within $\eps$ of the associated value. We will need to consider several cases separately.

\begin{case}All four values are the same side of $z^*$.\end{case}
Without loss of generality, say they are all at most $z^*$. Then taking $z_1,z_2$ to be the largest and smallest values, Lemma \ref{lem:nozero} gives $z_1>0$. Since $R(z)>z$ for all $z\in(0,z^*)$, Lemma \ref{lem:shift-up} gives a contradiction.

\begin{case}The highest and lowest values are on opposite sides of $z^*$ and correspond to different communities.\end{case}
Without loss of generality, we have $z_1^-\leq z_2^-$ and $z_1^+\leq z_2^+$ with $z_1^-<z^*<z_2^+$. Since clearly in this case $0<z^*<1$, we have $0<z_1^-<z^*<z_2^+<1$ by Lemma \ref{lem:nozero} (noting that non-negative determinant implies diagonal entries are positive).

For $i=1,2$, let\[
f_i(x;a,b)=\frac{\alpha_{1,i}xa+\alpha_{2,i}(1-x)b}{\alpha_{1,i}x+\alpha_{2,i}(1-x)},\]
so that $f_i(\hier[Y]t1;\hier[Z]t1,\hier[Z]t2)$ gives the probabilities of a neighbour of a new vertex added at time $t+1$ being red conditional on the new vertex being in community $i$. Note that $(f_1(x;a,b),f_2(x;a,b))$ is the outcome of multiplying $(a,b)$ by a matrix $\boldsymbol B_x$ obtained from $\boldsymbol A$ by multiplying first each row, and then each column, by positive constants. Since these operations preserve the sign of the determinant, it follows that $f_1(x;a,b)\leq f_2(x;a,b)$ whenever $a\leq b$. 

Recall that $\hier[Y]t1\to\nu_1$ almost surely, where $\nu_1\in(0,1)$. Now we have \[f_1(\nu_1;z_1^-,z_2^+)\leq f_2(\nu_1;z_1^-,z_2^+),\] and in particular either $f_1(\nu_1;z_1^-,z_2^+)\leq z^*$ or $f_2(\nu_1;z_1^-,z_2^+)\geq z^*$ (or both). Assume without loss of generality that $f_1(\nu_1;z_1^-,z_2^+)\leq z^*$. There is some $\eps>0$ such that $R(z)\geq z_1+\eps$ for $z\in[z_1^--\eps,z^*+\eps]$. We may choose $\delta>0$ such that $\abs{x-\nu_1}<\delta$ implies $\|B_{x}-B_{\nu_1}\|_\infty<\eps/4$, and consequently $\abs{f_1(x;a,b)-f_1(\nu_1;a,b)}\leq\eps/2$ for any $a,b$.

Suppose that $\hier[Z]t1\leq z_1+\eps/2$, $\hier[Z]t2\leq z_2+\eps/2$ and $\abs{\hier[Y]t1-\nu_1}<\delta$. Set $\hier[W]ti=f_i(\hier[Y]t1;\hier[Z]t1,\hier[Z]t2)$, so that the probability of a new vertex added to community $i$ at time $t$ being red is given by $R(\hier[W]ti)$. 

Suppose that $\hier[Z]t1\leq z_1^-+\eps/2$. Provided $t$ is sufficiently large that $\hier[Z]ti\in(z_1^--\eps/2,z_2^++\eps/2)$ for each $i$ and $|\hier[Y]t1-\nu_1|<\delta$, we have $\hier[W]t1\leq y'_1+\eps/2\leq z^*+\eps$, and also $\hier[W]t1\geq z_1^--\eps/2$, so consequently $R(\hier[W]t1)\geq z_1+\eps$. 

It follows that, for $t$ sufficiently large, while $\hier[Z]t1$ remains below $z_1^-+\eps/2$ the expected proportion of red edge-ends among new edge-ends meeting new vertices in community $1$ is at least $z_1^-+\eps$, and so the expected proportion of red edge-ends among new edge-ends in community $1$ is at least $\hier[Z]t1+\mu_1\eps/2$. Consequently, $\hier[Z]t1$ almost surely reaches $z_1+\eps/2$ from any point below it, but fails to reach $z_1+\eps/4$ from $z_1+\eps/2$ with positive probability. Therefore $\liminf_{t\to\infty}\hier[Z]t1\geq z_1+\eps/4$, a contradiction.

\begin{case}The highest and lowest values are on opposite sides of $z^*$ and correspond to the same community.\end{case}
Without loss of generality we have $z_1^-<z_2^-$ and $z_2^+<z_1^+$ with $0<z_1^-<z^*<z_1^+<1$ (where again Lemma \ref{lem:nozero} excludes the possibility of $0$ or $1$). Now we have
\[f_1(\nu_1;z_1^-,z_2^+)\leq f_2(\nu_1;z_1^-,z_2^+).\]
If $f_1(\nu_1;z_1^-,z_2^+)\leq z^*$, the analysis above goes through, so we must have 
\[z^*<f_1(\nu_1;z_1^-,z_2^+)\leq f_2(\nu_1;z_1^-,z_2^+)<z_2^+.\]
Similarly, we must have 
\[z_2^-<f_2(\nu_1;z_1^+,z_2^-)\leq f_1(\nu_1;z_1^+,z_2^-)<z^*.\]
In particular, it follows that $z_2^-<z^*<z_2^+$.

Write $y_1=f_1(\nu_1;z_1^-,z_2^+)$ and consider the minimum value of $R(z)$ for $z\in[z^*,y_1]$. If this is greater than $z_1^-$ then we have some $\eps>0$ such that provided $\hier[Z]t1\in (z_1^--\eps/2,z_1^-+\eps/2)$ and $\hier[Z]t2<z_2^++\eps/2$ and $|\hier[Y]t1-\nu_1|<\delta$, we have $\hier[W]t1\leq y_1+\eps$ and hence $R(\hier[W]t1)>z_1^-+\eps$, giving a contradiction as before.

Thus we must have some value $y\in(z^*,y_1)$ with $R(y)<z_2^-$. Choose some $y'\in(y,z_2^+)$ such that $f_2(\nu_1;z_1^-,y')>z^*$; this exists since $f_2(z_1^-,z_2^+)>z^*$. Take $\eta>0$ sufficiently small that $z^*,y,y',z_2^+$ are all at least $\eta$ apart, that $R(z)<z_2^--\eta$ for all $z\in[y-\eta,y+\eta]$, that $R(z)<y'-\eta$ for all $y'<z<y'+\eta$, and that $f_2(\nu_1;z^*,y+\eta)<y-\eta$.

Note that for any $\delta>0$ we have $\hier[Z]t2\geq z_1^--\delta$ and $\abs{\hier[Y]t1-\nu_1}<\delta$ for sufficiently large $t$, and $\hier[Z]t2\geq z_2^+-\delta$ infinitely often, and so choosing $\delta$ suitably we have $\hier[W]t2>y'+\eta$ infinitely often by choice of $y_1$. Similarly, since $f_2(\nu_1;z_1^+,z_2^-)<z^*$ we have $\hier[Z]t2,\hier[W]t2<y-\eta$ infinitely often. However, we will argue that the probability of reaching the region $\hier[z]t2>y'+\eta$ before returning to the region $\hier[Z]t2,\hier[W]t2<y-\eta$ is exponentially small, and thus almost surely this happens only finitely often, giving a contradiction.

Suppose $\hier[Z]{t_0}2,\hier[W]{t_0}2<y-\eta$, and consider $t\geq t_0$. We first argue that with suitably high probability we have $\hier[Z]t1<y$ before either $\hier[Z]t2>y$ or $\hier[W]t2>y$; by assumption one of these things eventually happens. Indeed, while $\hier[Z]t1>y$ and $\hier[W]t2<y$ we cannot have $\hier[Z]t2>y$ (since $\hier[W]t2$ is a convex combination of the other two). However, if $\hier[Z]t1>y$ and $\hier[Z]t2<y$ and $\hier[W]t2\in(y-\eta,y)$ then we have $\hier[Z]t1\geq\hier[W]t1\geq \hier[W]t2>y-\eta$, and hence $R(\hier[W]t2)<z_2^--\zeta<\hier[Z]t2$ and $R(\hier[W]t1)<\hier[Z]t1$. Thus both $\hier[Z]t1$ and $\hier[Z]t2$ are supermartingales, but in order for $\hier[W]t2$ to reach $y$ from $y-\eta/2$ before returning below $y-\eta$, one of them must increase by at least a constant, which is exponentially unlikely.

Now suppose $\hier[Z]{t_1}1,\hier[Z]{t_1}2,\hier[W]{t_1}2<y$. We next argue that we have at least one of $\hier[Z]t1<z^*$ or $\hier[Z]t2<z^*$ before either $\hier[Z]t1>y+\eta$ or $\hier[Z]t2>y+\eta$. Indeed, suppose $\hier[Z]t1,\hier[Z]t2\in[z^*,y+\eta]$. Then the same bounds apply to $\hier[W]t1,\hier[W]t2$, and so $R(\hier[W]t1),R(\hier[W]t2)\leq y-\eta$. Consequently the probability of either community reaching $y+\eta$ before one reaches $z^*$ is exponentially small.

If $\hier[Z]t2<z^*$ and $\hier[Z]t1<y+\eta$ then by choice of $\eta$ we have $\hier[Z]t2,\hier[W]t2<y-\eta$, as required. So we may assume $\hier[Z]{t_2}1<z^*$ and $\hier[Z]{t_2}2<y+\eta$. Now we argue that with suitably high probability we have $\hier[Z]t2,\hier[W]t2<z^*$ before $\hier[Z]t2>y'+\eta$. First note that if $\hier[Z]t2\in(y',y'+\eta/2)$ and $\hier[Z]t1<\max\{\hier[Z]t2,y'\}$ we have $z^*<\hier[W]t2<y'+\eta/2$ and thus $R(\hier[W]t2)<y'$, so with high probability we have $\hier[Z]t1>y'$ before $\hier[Z]t2>y'+\eta$. However, if $z^*<\hier[W]t1<y'+\eta$ then we have $R(\hier[W]t1)<y'$, and so it is exponentially unlikely for $\hier[Z]t2$ to reach $y'+\eta$ unless either $\hier[Z]t2$ exceeds $y'+\eta$ (which we know does not happen first) or $\hier[W]t1<z^*$ while $\hier[Z]t1>z^*$. However, these latter conditions imply $\hier[Z]t2,\hier[W]t2<z^*$, as required.

\subsubsection{Proof of Theorem \ref{same2}}\label{period2point}

Recall that here $A=\begin{pmatrix} 0 & 1 \\ 1 & 0\end{pmatrix}$, $z^{*}$ is the only stable fixed point of $R$, and we have $\tilde{z}\neq z^*$, a linearly stable fixed point of the second iterate of $R$, that is such that $R^2(\tilde{z})=R(R(\tilde{z}))=\tilde{z}$ and that $(R^2)'(\tilde{z})<1$. See Figure \ref{fig:rcond} for an example. We will show that there is positive probability that $\hier[Z]t1\to \tilde{z}$ and $\hier[Z]t2\to R(\tilde{z})$; since $\tilde z\neq z^*$ these are not equal.

We start by calculating the form of the vector field $F$ for this community structure.  Defining, for $i=1,2$ and $j=1,2$, $Z_{i,j}(z)=\frac{x_{i,j}}{x_{i,1}+x_{i,2}}$, for $j=1,2$ we have \[Q_{1,j}(x)=\frac{x_{2,j}}{x_{2,1}+x_{2,2}}=Z_{2,j}(x),\] and similarly $Q_{2,j}(x)=Z_{1,j}(x)$.  Hence $G_{i,1}(x)=m\mu_i R(Z_{3-i,1}(x))$ and, similarly $G_{i,2}(x)=m\mu_i (1-R(Z_{3-i,2}(x)))$.

We also have, for $i,j=1,2$, \[H_{i,j}(x)=m\mu_{3-i} \frac{x_{i,j}}{x_{i,1}+x_{i,2}}=m\mu_{3-i} Z_{i,j}(x),\] and so we have \[F_{i,1}(x)=m \left(\mu_i R(Z_{3-i,1}(x))+\mu_{3-i} Z_{i,1}(x)\right)-2mx_{i,1}\] and \[F_{i,2}(x)=m\left(\mu_i(1-R(Z_{3-i,1}(x)))+\mu_{3-i}Z_{i,2}(x)\right)-2mx_{i,2}.\]

Let $x=\frac12(\tilde{z},1-\tilde{z},R(\tilde{z}),1-R(\tilde{z}))$.  Then $x$ is a fixed point of $F$; indeed $Z_{1,1}(x)=\tilde{z}$ and $Z_{2,1}(x)=R(\tilde{z})$, and we can check that $F_{1,1}(z)=m(\mu_i R^2(\tilde{z})+\mu_{3-i}\tilde{z}-2\tilde{z})=0$ by our assumption on $\tilde{z}$ and because $\mu_i+\mu_{3-i}=1$, and similarly for the other components of $F$.

We now investigate the stability of $x$ as a fixed point of $F$.  For this community structure we can see that the equations \eqref{edgeendlim} giving the measure $\nu$ reduce to $\nu_1=\nu_2=\frac12$, whatever the values of $\mu_1$ and $\mu_2$.  Hence, using Lemma \ref{lem:nu}, both $\hier t{1,1}+\hier t{1,2}$ and $\hier t{2,1}+\hier t{2,2}$ converge almost surely to $\frac12$ regardless of the types in the two communities, and so we may consider the dynamics of $F$ restricted to the set $\{x\in\Delta^3: x_{1,1}+x_{1,2}=x_{2,1}+x_{2,2}=\frac12\}$.  On this set we have $Z_{i,j}(x)=2x_{i,j}$ and so $\frac{\partial F_{i,1}}{\partial x_{i,1}}=-2\mu_1m$ and $\frac{\partial F_{i,1}}{\partial x_{3-i,1}}=2\mu_i mR'(2x_{2,1})$.  Hence for the choice of $x$ above the Jacobian is \[m\begin{pmatrix}-2\mu_1 & 2\mu_1 R'(R(\tilde{z})) \\ 2\mu_2 R'(\tilde{z}) & -2\mu_2 \end{pmatrix},\] which has both eigenvalues negative if and only if $(R^2)'(\tilde{z})=R'(\tilde{z})R'(R(\tilde{z}))<1$.  Hence if $(R^2)'(\tilde{z})<1$ then $x$ is a linearly stable fixed point of $F$.

The conclusion follows by applying Theorem 2.16 of Pemantle \cite{pemantlesurvey}, observing that there will always be a sequence of community and type for each vertex which has positive probability and gets arbitrarily close to $x$.

\subsubsection{Proof of Theorem \ref{determinant}, part 2}

Now suppose we have a matrix for which $\det(A)<0$. In particular this implies $B_{\boldsymbol \nu}\binom10$ has second coordinate larger. We can choose $\eps>0$ and $z_1<z_2$ such that given $\mathbf{x}$ is within $\eps$ of $\binom10$ and $\boldsymbol{\nu}'$ is within $\eps$ of $\boldsymbol{\nu}$ then $\mathbf{y}:=B_{\boldsymbol\nu'}\mathbf x$ satisfies $y_1<z_1<z_2<y_2$. We can choose $r,m$ such that $z_1<r/m<z_2$ and $\prob{\bin{m,z_1}>r},\prob{\bin{m,z_2}\leq r}<\eps$. We then use a modified minority rule: $p_k=1$ if $k\leq r$ and $p_k=0$ otherwise. Note that, since $p_k$ is decreasing, so is $R(z)$, meaning that $P'(z)\leq-1$, and so necessarily $P$ has a single root which is linearly stable.

By choice of parameters, provided $\hier[Z]t1\geq 1-\eps$, $\hier[Z]t2<\eps$, and $|\hier[Y]t1-\nu_1|<\eps$ then this trend in $\hier[Z]ti$ is preserved and with positive probability stays true forever, so they don't converge to the same value.

\subsection{Proof of Proposition \ref{threecycle}}

Recall that here we have $S=\{1,2,3\}$, \[A=\begin{pmatrix} 0 & 1 & 0 \\ 0 & 0 & 1 \\ 1 & 0 & 0\end{pmatrix}\] and $\mu_1=\mu_2=\mu_3=\frac13$.  We start off by identifying the vector field $F$ in this community structure; this calculation is similar to that in Section \ref{period2point}.  We have, for $j=1,2$, \[Q_{1,j} = \frac{x_{3,j}}{x_{3,1}+x_{3,2}} = Z_{3,j}(x),\] where for $i=1,2,3$ and $j=1,2$ we define $Z_{i,j}(x)=\frac{x_{i,j}}{x_{i,1}+x_{i,2}}$, and similarly $Q_{2,j}=Z_{1,j}(x)$ and $Q_{3,j}=Z_{2,j}(x)$.  Hence $G_{1,1}(x)=\frac{m}{3}R(Z_{3,1}(x))$, $G_{2,1}(x)=\frac{m}{3}R(Z_{1,1}(x))$ and $G_{3,1}(x)=\frac{m}{3}R(Z_{2,1}(x))$.   We also have, for $i=1,2,3$ and $j=1,2$, $H_{i,j}(x)=\frac{m}{3}\frac{x_{i,j}}{x_{i,1}+x_{i,2}}$, again similarly to Section \ref{period2point}.  Thus 
\begin{align*} F_{1,1}(x) &= \frac{m}{3}(R(Z_{3,1}(x))+Z_{1,1}(x))-2mx_{1,1} \\
F_{2,1}(x) &= \frac{m}{3}(R(Z_{1,1}(x))+Z_{2,1}(x))-2mx_{2,1} \\
F_{3,1}(x) &= \frac{m}{3}(R(Z_{2,1}(x))+Z_{3,1}(x))-2mx_{3,1},\end{align*} with similar expressions for the $F_{i,2}(x)$.

Next, we identify the fixed points of $F$.  Similarly to Section \ref{period2point}, symmetry implies that the measure $\nu$ from \cite{geopref9} has $\nu_1=\nu_2=\nu_3=\frac13$, so again Lemma \ref{lem:nu} means that we can consider the dynamics of $F$ restricted to the set $\Delta_{\nu}=\{x\in \Delta^5: x_{1,1}+x_{1,2}=x_{2,1}+x_{2,2}=x_{3,1}+x_{3,2}=\frac13\}$, meaning that we can take $Z_{i,j}(x)=3x_{i,j}$.  It is then easy to check that any $x\in\Delta_{\nu}$ for which $F(x)=0$ must satisfy $3x_{1,1}=R^3(3x_{1,1})$, and because for our choice of rule $R$ is decreasing the only fixed point is $x=(\frac{1}{6},\frac{1}{6},\frac{1}{6},\frac{1}{6},\frac{1}{6},\frac{1}{6})$.

To investigate whether this fixed point is a possible limit, we start by checking its stability.  At this value of $x$ the Jacobian of $F$ restricted to $\Delta_{\nu}$ is 
\[m\begin{pmatrix}-1 & 0 & R'(\frac12) \\ R'(\frac12) & -1 & 0 \\ 0 & R'(\frac12) & -1\end{pmatrix},\]
whose eigenvalues are $-1$ and $-1-\frac12R'(\frac12)\pm i\frac{\sqrt{3}}{2}R'(\frac12)$.  It follows that this stationary point of $F$ is linearly unstable if $R'(\frac12)<-2$.  For the minority rule with odd $m$, $R(z)=\sum_{j=0}^{(m-1)/2}\binom{m}{j}z^j(1-z)^{m-j}$, with $R'(z)=-m\binom{m-1}{(m-1)/2}z^{\frac{m-1}{2}}(1-z)^{\frac{m-1}{2}}$, so $R'(\frac12)=\frac{-m!}{2^{m-1}((m-1)/2)!)^2}$, which is decreasing in $m$ and is less than $-2$ for $m\geq7$, so in these cases the fixed point in linearly unstable.

In those cases where $(\frac{1}{6},\frac{1}{6},\frac{1}{6},\frac{1}{6},\frac{1}{6},\frac{1}{6})$ is the only stationary point of $F$ and is linearly unstable, we can apply Theorem 9.1 of Bena\"{\i}m \cite{benaim} to show that convergence to the fixed point $(\frac{1}{6},\frac{1}{6},\frac{1}{6},\frac{1}{6},\frac{1}{6},\frac{1}{6})$ has probability zero.  To apply this result we need to show that within a neighbourhood of the fixed point the expectation of the positive part of the noise term $\xi^{(t)}$ in any direction is uniformly bounded away from zero by a constant.  To note that this holds for the minority rule, observe that for $\mathbf{X}_t$ in a suitable neighbourhood of $(\frac{1}{6},\frac{1}{6},\frac{1}{6},\frac{1}{6},\frac{1}{6},\frac{1}{6})$ the probability that the new vertex is in community $i$ is $\frac13$ and that conditional on this the probability of it being type $j$ is at least $\frac12-\eps$, while if both these events happen the increase in the number of edge ends at vertices of type $j$ in community $i$ is at least $m$, compared with an expected increase of at most $2m\left(\frac16+\eps\right)$.  Therefore $\EE((\xi_{i,j}^{t+1})^+|\salj_t)\geq \frac{m}{12}$, say, for a suitably small choice of $\eps$, and as this applies to all choices of $i$ and $j$ the same bound will apply to other directions.  This gives the stated result for the minority rule with $m\geq 7$.

\subsection{Proof of Theorem \ref{linear}}

Define $\hier[\hat{Z}]ti=\hier t{i,1}/\nu_i=\hier[Z]ti(\hier[Y]ti/\nu_i)$, where $\nu_i$ is the almost sure limit of $\hier[Y]ti$ as $t\to\infty$ from Proposition 3.1 of \cite{geopref9}; note that because of this convergence it is sufficient to show convergence of the $\hier[\hat{Z}]ti$.We will start by giving a result on the rate of convergence of the $\hier[Y]ti$ to $\nu_i$, applying general results on the rates of convergence of Robbins-Monro stochastic approximation processes to the process in \cite{geopref9} which gives the values of the $\hier[Y]ti$.

\begin{lemma}\label{lem:ratecon}
There exists $r>0$ such that \[\sum_{i=1}^N \left|\hier[Y]ti-\nu_i\right|=O(t^{-r})\text{ as }t\to\infty.\]
\end{lemma}

\begin{proof}
We use Proposition 2.3 of Dereich and M\"{u}ller-Gronbach \cite{dereich2019general}, for which we need to check that a number of conditions are satisfied, labelled (I) to (III) and A.1 and A.2 in \cite{dereich2019general}.  Conditions (I) to (III) follow from the fact that we have a stochastic approximation process, and in particular, following their notation except for denoting time by $t$, we have $R_t=0$ and $\gamma_t=\frac{1}{t+1+e_0/m}$, where $e_0$ is the number of edges in the initial graph.  That $R_t=0$ means we can also take the $\eps_t$ in the notation of \cite{dereich2019general} (and hence also $e_{k,t}(r)$) as being zero.  That condition A.2 is satisfied follows from the boundedness of the noise terms.  That condition A.1 is satisfied with some $L>0$ follows from the strict convexity of the Lyapunov function in \cite{geopref9}.
	
To apply Proposition 2.3 of \cite{dereich2019general} we need to estimate $\tau_{k,t}(r)=\prod_{j=k+1}^t (1-\gamma_j r)=\Theta((t/k)^{-r})$ and $s_{k,t}^2(r)=\sum_{j=k}^t \gamma_j ^2\tau_{j,t}(r)^2=\Theta(t^{2r}k^{-2r-1})$.  Then for $0<r<\min\{\frac12,L\}$, $\tau_{j,t}(r)=\Theta((t/j)^{-r})$ and $s_{k,t}^2(r)=\Theta(t^{-2r})$ so Proposition 2.3 of \cite{dereich2019general} thus gives that the $L^p$ norm of the difference between the stochastic approximation process and its limit is $O(t^{-r})$.
\end{proof}

Next, we show that there exists a linear combination of the $\hier[\hat{Z}]ti$ which converges to a random limit.

\begin{lemma}\label{lem:ctsmarkov}There exist $\sigma_i, 1\leq i\leq N$ with all $\sigma_i>0$ and $\sum_{i=1}^N \sigma_i=1$ and a random variable $M$ such that $\sum_{i=1}^N \sigma_i \hier[\hat{Z}]ti \to M$, almost surely, as $t\to\infty$.\end{lemma}

\begin{proof}In the linear model, for any community $i$ the expected number of new red edge ends in that community at existing vertices is $m\sum_{j=1}^N\mu_j \frac{\alpha_{i,j}\hier[Y]ti\hier[Z]ti}{\sum_{k=1}^N\alpha_{k,j}\hier[Y]tk}$, while the expected number of new red edge ends in the community at the new vertex is $m\mu_i \frac{\sum_{j=1}^N \alpha_{j,i}\hier[Y]tj \hier[Z]tj}{\sum_{j=1}^N \alpha_{j,i} \hier[Y]tj}$.    It follows that 
\begin{multline}\label{zhat}\EE(\hier[\hat{Z}]{t+1}i\mid\salj_t)-\hier[\hat{Z}]ti=\\
\frac{m}{2m(t+n_0+1)\nu_i}\Biggl(\mu_i \frac{\sum_{i=1}^N \alpha_{j,i}\hier[Y]tj \hier[Z]tj}{\sum_{i=1}^N \alpha_{j,i} \hier[Y]tj}+\sum_{j=1}^N\mu_j \frac{\alpha_{i,j}\hier[Y]ti\hier[Z]ti}{\sum_{k=1}^N\alpha_{k,j}\hier[Y]tk}-2\nu_i\Biggr).\end{multline}
	
From the form of the function $G$ driving the stochastic approximation in \cite{geopref9} and that $G(\nu)$ is zero for the limiting measure $\nu$, it follows that \[\sum_{j=1}^N\mu_j \frac{\alpha_{i,j}\nu_i}{\sum_{k=1}^N\alpha_{k,j}\nu_k}=2\nu_i-\mu_i.\]  Therefore we can write the term in brackets on the right hand side of \eqref{zhat} as 
\[\mu_i \frac{\sum_{i=1}^N \alpha_{j,i}\hier[Y]tj \hier[Z]tj}{\sum_{i=1}^N \alpha_{j,i} \hier[Y]tj}-\mu_i\hier[Z]ti+\phi_i(Y^{(t)},\nu)\hier[Z]ti,\]
where for two probability measures $\xi$ and $\eta$ on $S$ \[
\phi_i(\xi,\eta)=\sum_{j=1}^N\mu_j\left(\frac{\alpha_{i,j}\xi_i}{\sum_{k=1}^N\alpha_{k,j}\xi_k}- \frac{\alpha_{i,j}\eta_i}{\sum_{k=1}^N\alpha_{k,j}\eta_k}\right).\]
	
Hence for any linear combination of the communities we have \begin{multline*}\EE\left(\sum_{i=1}^N s_i \hier[\hat{Z}]{t+1}i\biggm|\salj_t\right)-\sum_{i=1}^N s_i\hier[\hat{Z}]ti=
\sum_{i=1}^N s_i\phi_i(Y^{(t)},\nu)\hier[Z]ti\\
+\frac{m}{2m(t+n_0+1)}\sum_{i=1}^N s_i\mu_i
\left(\frac{\sum_{j=1}^N \alpha_{j,i}\hier[Y]tj \hier[Z]tj}{\sum_{j=1}^N \alpha_{j,i} \hier[Y]tj}-\hier[Z]ti\right).\end{multline*}
	
We aim to find a particular linear combination with coefficients $\sigma_i$, $i=1,\ldots, N$, such that $\sum_{i=1}^N \sigma_i \hier[\hat{Z}]ti$ would be a martingale if $Y^{(t)}$ were equal to its limit $\nu$.  Because we would then have $\hier[\hat{Z}]tj=\hier[Z]tj$ for all $j$, this will be the case if for any values of $x_1,x_2,\ldots,x_N$ we have \[\sum_{i=1}^N \sigma_i \mu_i\left(\frac{\sum_{j=1}^N \alpha_{j,i} \nu_j x_j}{\sum_{j=1}^N \alpha_{j,i} \nu_j}-x_i\right)=0,\] which is equivalent to \[\sum_{i=1}^N \sigma_i\frac{\mu_i}{\sum_{j=1}^N \alpha_{j,i} \nu_j}\left(\sum_{j=1}^N \alpha_{j,i} \nu_j(x_j-x_i)\right)=0.\]  Define an $N\times N$ matrix $\Psi$ by letting its entries $\psi_{i,j}=\frac{\mu_i \alpha_{j,i}\nu_j}{\sum_{k=1}^N \alpha_{k,i}\nu_k}$ if $i\neq j$ and $\psi_{i,i}=\frac{\mu_i \alpha_{i,i}\nu_i}{\sum_{k=1}^N \alpha_{k,i}\nu_k}-\mu_i$.  Then $\Psi$ is the generator matrix of a continuous time Markov chain on the space of communities with transition rates from community $i$ to community $j$ given by the strength of the connection between the two communities times the proportion of vertices which are in community $i$ times the proportion of edge ends which are in community $j$.  Under our assumptions on the community structure this chain has a unique stationary distribution, which we write as $(\sigma_1,\sigma_2,\ldots,\sigma_N)$.  It follows immediately that $\sum_{i=1}^N \sigma_i\psi_{i,j}=0$ for all $j$, thus satisfying the condition above.
	
Write $M_t=\sum_{i=1}^N \sigma_i \hier[\hat{Z}]ti$.  Then, choosing suitable values of $\eps$ and $K$ and assuming $\hier[Y]ti>\eps$, which occurs for sufficiently large times if $\eps$ is small enough, \[\bigl|\EE(M_{t+1}\mid\salj_t)-M_t\bigr|\leq Kt^{-1}\sum_{i=1}^N \left|\hier[Y]ti-\nu_i\right|.\]  Define $\tilde{M}_t=M_t-\sum_{j=2}^m (\EE(M_{j}\mid\salj_t)-M_{j-1})$; it immediately follows from the above that $(\tilde{M}_t)$ is a martingale, and by Lemma \ref{lem:ratecon} the expectation of the difference between the increments of $(M_t)$ and $(\tilde{M}_t)$ at time $t$ is $O(t^{-(1+r)})$, which is summable.  Hence $M_t$ converges almost surely, to a random limit $M$, as $t\to\infty$.\end{proof}

To complete the proof of Theorem \ref{linear}, we need to show that the proportions of red in each community all converge to the same limit $M$.  To do this, note that, because $\Psi$ from the proof of Lemma \ref{lem:ctsmarkov} is the generator matrix of a finite state continuous time Markov chain, for any row vector $s\in \RR^N$ we have that the sum of the elements in $s\Psi$ is zero, and thus the $(N-1)$-dimensional subspace of row vectors whose elements sum to zero is preserved under $\Psi$.  Furthermore, by our assumptions on the community structure, the multiplicity of $0$  as an eigenvalue of $\Psi$ is $1$ and all other eigenvalues are strictly negative; let $-\lambda_2<0$ be the largest eigenvalue of $\Psi$ other than $0$ and let $-\lambda_N$ be the smallest eigenvalue of $\Psi$.

Then for any row vector $s\in \RR^N$ whose elements sum to zero, we have 
\[\EE\left(\sum_{i=1}^N s_i\hier[\hat{Z}]{t+1}i\Biggm|\salj_t\right)-\sum_{i=1}^N s_i\hier[\hat{Z}]ti \leq \frac{-\lambda_2 m}{2m(t+n_0+1)}\sum_{i=1}^N s_i \hier[\hat{Z}]ti+Kt^{-1}\sum_{i=1}^N \left|\hier[Y]ti-\nu_i\right|,\] and \[\EE\left(\sum_{i=1}^N s_i\hier[\hat{Z}]{t+1}i\Biggm|\salj_t\right)-\sum_{i=1}^N s_i\hier[\hat{Z}]ti \geq \frac{-\lambda_N m}{2m(t+n_0+1)}\sum_{i=1}^N s_i \hier[Z]ti-Kt^{-1}\sum_{i=1}^N \left|\hier[Y]ti-\nu_i\right|,\] again for a suitable choice of $K$ and sufficiently large times.  By using inequality versions of standard stochastic approximation results (for example Lemma 5.4 of Jordan and Wade \cite{jordanwade}) it follows that $\sum_{i=1}^N s_i\hier[\hat{Z}]ti \to 0$, almost surely.  For $i\neq j$, this applies to $s=e_i-e_j$ (where $e_i$ and $e_j$ are standard unit vectors) and so shows that $\hier[Z]ti-\hier[Z]tj\to 0$ almost surely.  Putting this together with our previous conclusion on $\sum_{i=1}^N \sigma_i\hier[\hat{Z}]ti$, we have $\hier[Z]ti \to M$ almost surely as $t\to\infty$, which completes the proof of Theorem \ref{linear}.

\section{Examples}\label{examples}

We now consider some examples of community structures and type assignment rules where we can apply our results and compare them with the results for the single community model which can be deduced from \cite{AMR}; we also show some simulations which illustrate them.  In some cases we can calculate stationary points explicitly and hence give explicit phase transitions. Figure \ref{fig:four-plots} shows the function $R(z)$ for some examples we consider.

    \begin{figure}
        \centering
        \begin{subfigure}[b]{0.475\textwidth}
            \centering
	\begin{tikzpicture}
		\begin{axis}[
			width=.9\textwidth,
			height=4cm,
			axis lines = left,
			ymin=0,ymax=1,
			xlabel = $z$,
			xmin=0, xmax=1,
			ylabel = $R(z)$,
			]
			\addplot [
			domain=0:1, 
			samples=300, 
			color=black,
			thick,
			]
			{-2*x^3 + 3*x^2};
			\addplot [
			domain=0:1, 
			samples=300, 
			color=black,
			]
			{x};						
		\end{axis}
	\end{tikzpicture}
	\caption{Majority of three (Section \ref{majwins}).}
            \label{fig:4p1}
        \end{subfigure}%
        \hfill%
        \begin{subfigure}[b]{0.475\textwidth}  
            \centering 
	\begin{tikzpicture}
		\begin{axis}[
			width=.9\textwidth,
			height=4cm,
			axis lines = left,
			ymin=0,ymax=1,
			xlabel = $z$,
			xmin=0, xmax=1,
			ylabel = $R(z)$,
			]
			\addplot [
			domain=0:1, 
			samples=300, 
			color=black,
			thick,
			]
			{x^3 - 1.5*x^2 + 1.5*x};
			\addplot [
			domain=0:1, 
			samples=300, 
			color=black,
			]
			{x};						
		\end{axis}
	\end{tikzpicture}
	\caption{Random visible type (Section \ref{randvis}).}
            \label{fig:4p2}
        \end{subfigure}%
        \vskip\baselineskip%
        \begin{subfigure}[b]{0.475\textwidth}   
            \centering 
	\begin{tikzpicture}
		\begin{axis}[
			width=.9\textwidth,
			height=4cm,
			axis lines = left,
			ymin=0,ymax=1,
			xlabel = $z$,
			xmin=0, xmax=1,
			ylabel = $R(z)$,
			]
			\addplot [
			domain=0:1, 
			samples=300, 
			color=black,
			thick,
			]
			{2*x^3 - 3*x^2+1};
			\addplot [
			domain=0:1, 
			samples=300, 
			color=black,
			]
			{x};						
		\end{axis}
	\end{tikzpicture}
		\caption{Minority of three (Section \ref{minority}).}
            \label{fig:4p3}
        \end{subfigure}%
        \hfill%
        \begin{subfigure}[b]{0.475\textwidth}   
            \centering 
	\begin{tikzpicture}
		\begin{axis}[
			width=.9\textwidth,
			height=4cm,
			axis lines = left,
			ymin=0,ymax=1,
			xlabel = $z$,
			xmin=0, xmax=1,
			ylabel = $R(z)$,
			]
			\addplot [
			domain=0:1, 
			samples=300, 
			color=black,
			thick,
			]
			{-2.25*x^3 + 3.75*x^2-0.75*x+0.25};
			\addplot [
			domain=0:1, 
			samples=300, 
			color=black,
			]
			{x};						
		\end{axis}
	\end{tikzpicture}
		\caption{Rule with a touchpoint (Section \ref{tp-rule}).}
            \label{fig:4p4}
        \end{subfigure}%
        \caption{Plots of $R(z)$ for the examples given in Sections \ref{majwins}, \ref{randvis}, \ref{minority} and \ref{tp-rule}.} 
        \label{fig:four-plots}
    \end{figure}

\subsection{Symmetric community structure, majority wins rule}\label{majwins}

One example of a community structure which is analysed in detail in \cite{geopref9} is where there is a two-point space $S=\{1,2\}$ with a simple symmetric attractiveness function given by, in the notation of Theorem \ref{different}, \[A_1=\begin{pmatrix}1 & 1 \\ 1 & 1 \end{pmatrix},\] leading to the matrix $A_{\theta}$ having the form \[A_{\theta}=\begin{pmatrix}1 & \theta \\ \theta & 1\end{pmatrix}.\] Here $\theta$ can be thought of as describing the strength of the connection between the two communities; in particular if $\theta=0$ the communities evolve separately, with no connections developing, while if $\theta=1$ vertices have no preference to connect to their own community or the other one so the graph effectively evolves as the single community model.  Here $\mu_1$ and $\mu_2=1-\mu_1$ give the relative sizes of the two communities.  In this case the limiting proportions $\nu_i$ of edge ends in each community can be found by solving (2.3) in \cite{geopref9}; there is no closed form solution in general, but in the case $\mu_1=\mu_2=1/2$ it is easy to see that $\nu_1=\nu_2=1/2$.

The \emph{majority wins} rule means that the new vertex connects to $m$ existing vertices and takes the type of the majority of them.  (If $m$ is even, then ties will be broken randomly.)  Here we consider it with $m=3$, but it can also be considered with other values of $m$.  This gives $p_0=p_1=0, p_2=p_3=1$, and so $R(z)=3z^2(1-z)+z^3$ with stable fixed points at $z=0$ and $z=1$ and an unstable fixed point at $z=\frac12$, see Figure \ref{fig:4p1}. Hence the results of \cite{AMR} show that for the single community model the proportion of red vertices tends to either 0 or 1, each with positive probability.  The question then arises of whether, with two communities, the proportion of red vertices in those two communities will tend to the same limit, and Theorem \ref{different} shows that there exists $\theta_{\mathrm{crit}}>0$ such that for $\theta<\theta_{\mathrm{crit}}$ there is positive probability that this does not happen and that different limits occur in the two communities.

In the case $\mu_1=\mu_2=\frac12$, a limit for the stochastic process $\mathbf{X}^{(t)}$ must be of the form $(y_1,\frac12-y_1,y_2,\frac12-y_2)$.  In this case, it is possible to explicitly solve the equations to find all the possible stationary points of this form and classify their stability; this work appears in the third author's PhD thesis \cite{yarrowthesis}.  The stationary points are given in Table \ref{statpoints}, where \begin{align*}
	S&=3\theta^3-9\theta^2+3\theta-1,\\
	R&=(\theta-1)\sqrt{-(7\theta-1)(\theta+1)^3},\\
	U&=(\theta-1)^3.
\end{align*}
\begin{table}[ht]
	\centering
	\begin{tabular}{lll}
		&$y_1$ & $y_2$\\\hline
		1)&$\frac{1}{4}+\frac{(\theta+1)}{4(\theta-1)}\sqrt{\frac{5\theta-1}{\theta-1}}$ & $\frac{1}{4}-\frac{(\theta+1)}{4(\theta-1)}\sqrt{\frac{5\theta-1}{\theta-1}}$\\
		2)&$\frac{1}{4}-\frac{(\theta+1)}{4(\theta-1)}\sqrt{\frac{5\theta-1}{\theta-1}}$ & $\frac{1}{4}+\frac{(\theta+1)}{4(\theta-1)}\sqrt{\frac{5\theta-1}{\theta-1}}$\\
		3)&$0$ & $0$\\
		4)&$\frac{1}{2}$ & $\frac{1}{2}$\\
		5)&$\frac{1}{4}$ & $\frac{1}{4}$\\
		6)&$\frac{1}{4}+\frac{\sqrt{2}}{8}\sqrt{\frac{S+R}{U}}$ & $\frac{1}{4}-\frac{\sqrt{2}}{8}\sqrt{\frac{S-R}{U}}$\\
		7)&$\frac{1}{4}-\frac{\sqrt{2}}{8}\sqrt{\frac{S+R}{U}}$ & $\frac{1}{4}+\frac{\sqrt{2}}{8}\sqrt{\frac{S-R}{U}}$\\
		8)&$\frac{1}{4}+\frac{\sqrt{2}}{8}\sqrt{\frac{S-R}{U}}$ & $\frac{1}{4}-\frac{\sqrt{2}}{8}\sqrt{\frac{S+R}{U}}$\\
		9)&$\frac{1}{4}-\frac{\sqrt{2}}{8}\sqrt{\frac{S-R}{U}}$ & $\frac{1}{4}+\frac{\sqrt{2}}{8}\sqrt{\frac{S+R}{U}}$
	\end{tabular}
	\caption{Stationary points for the flow for the symmetric $m=3$ majority wins model.}
	\label{statpoints}
\end{table}

The stationary points given by 3), 4) and 5) exist for all $\theta$.  Stationary points 1) and 2) are real and have $0\leq y_1\leq \frac12$ and $0\leq y_2\leq \frac12$ if and only if $0\leq \theta\leq \frac15$, while stationary points 6) to 9) are real and have $0\leq y_1\leq \frac12$ and $0\leq y_2\leq \frac12$ if and only if $0\leq \theta\leq \frac17$.

Restricted to the plane of points of the form $(y_1,\frac12-y_1,y_2,\frac12-y_2)$, the eigenvalues of the Jacobian of $F$ at $x=(x_{1,1},x_{1,2},x_{2,1},x_{2,2})=(y_1,\frac12-y_1,y_2,\frac12-y_2)$ are
\[
	\lambda=-3+\frac{9}{1+\theta}\left(J\pm\sqrt{J^2-4K}\right)
\]
where
\[
J=Q_{1,1}(x)\left(1-Q_{1,1}(x)\right)+Q_{2,1}(x)\left(1-Q_{2,1}(x)\right)
\]
and
\[
K=Q_{1,1}(x)Q_{2,1}(x)\left(1-Q_{1,1}(x)\right)\left(1-Q_{2,1}(x)\right)\left(1-\theta^2\right).
\]

Using this, the stationary points 3) and 4) are stable for all values of $\theta$.  Convergence to 3) represents type 2 dominating in both communities, and convergence to 4) represents type 1 dominating in both communities.  Meanwhile the stationary point 5), which corresponds to equal proportions of each type in both communities, is linearly unstable for all values of $\theta$, and stationary points 6) to 9) are linearly unstable for all values of $\theta$ for which they are meaningful.  However, the stationary points 1) and 2) turn out to be linearly stable when $\theta<\frac17$, and linearly unstable for $\frac17<\theta\leq \frac15$. 

Applying Theorem 2.16 of Pemantle \cite{pemantlesurvey} shows that the stable stationary points have positive probability of being limits and applying Theorem 2.17 of \cite{pemantlesurvey} shows that the linearly unstable stationary points have probability zero of being limits. Thus the stationary points 1) and 2) are limits with positive probability when $\theta<\frac17$ but with probability zero when $\theta>\frac17$.  These show the majority of edge ends in one community being of vertices of one type, and the majority of edge ends in the other community being of the other type; in each case there are a small (but $\Theta(t)$) number of vertices of the other community's type present.  Hence for this model the value $\theta_{\mathrm{crit}}$ in Theorem \ref{different} is $\frac17$: for $\theta<\frac17$ it is possible to have different limits in the two communities. 

\begin{figure}[ht]
	\centering\begin{tabular}{cc}
		\scalebox{0.4}{\includegraphics[trim=4cm 0.5cm 4cm 0.5cm, clip=true]{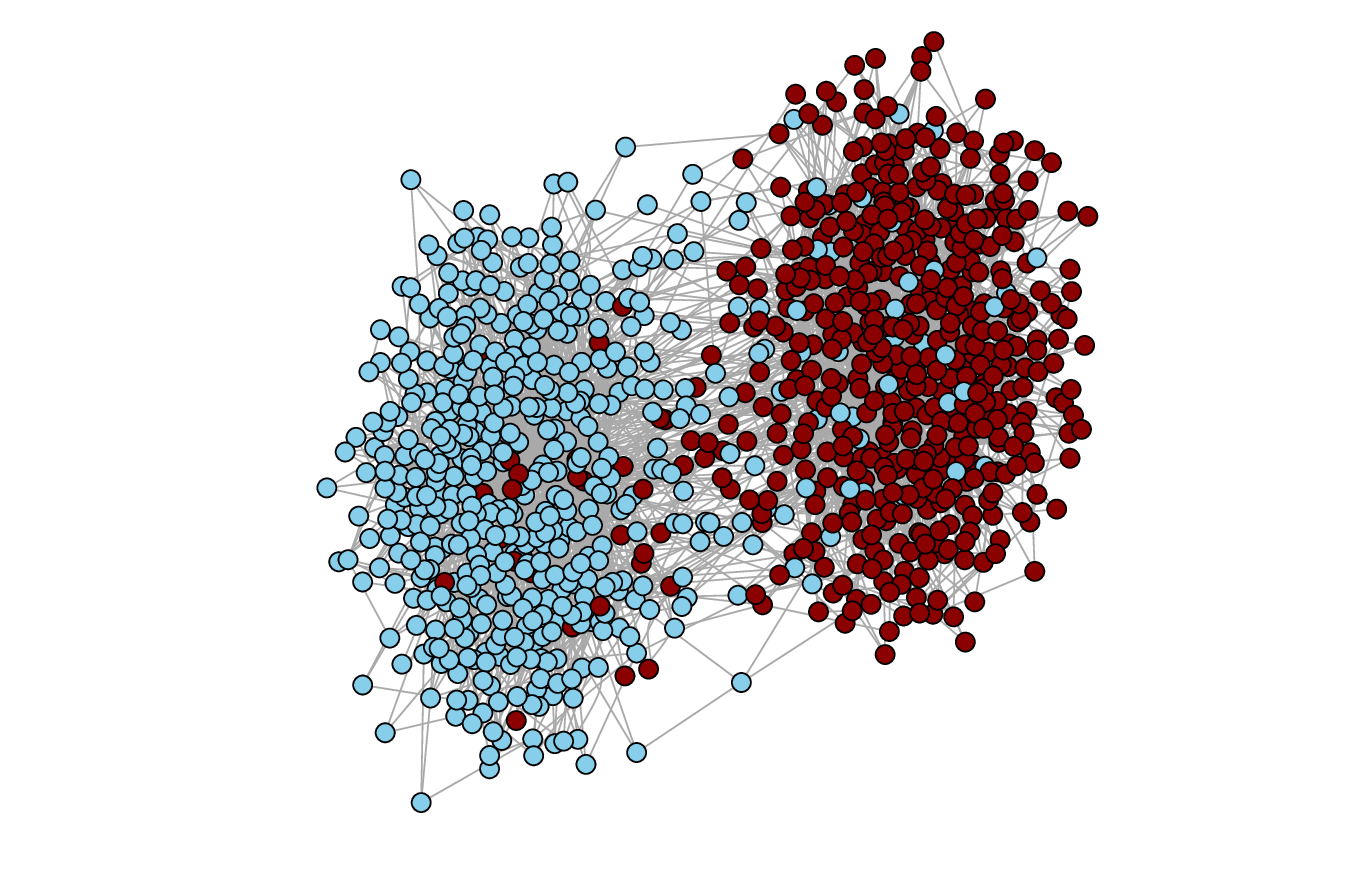}} & \scalebox{0.4}{\includegraphics[trim=4cm 0.5cm 4cm 0.5cm, clip=true]{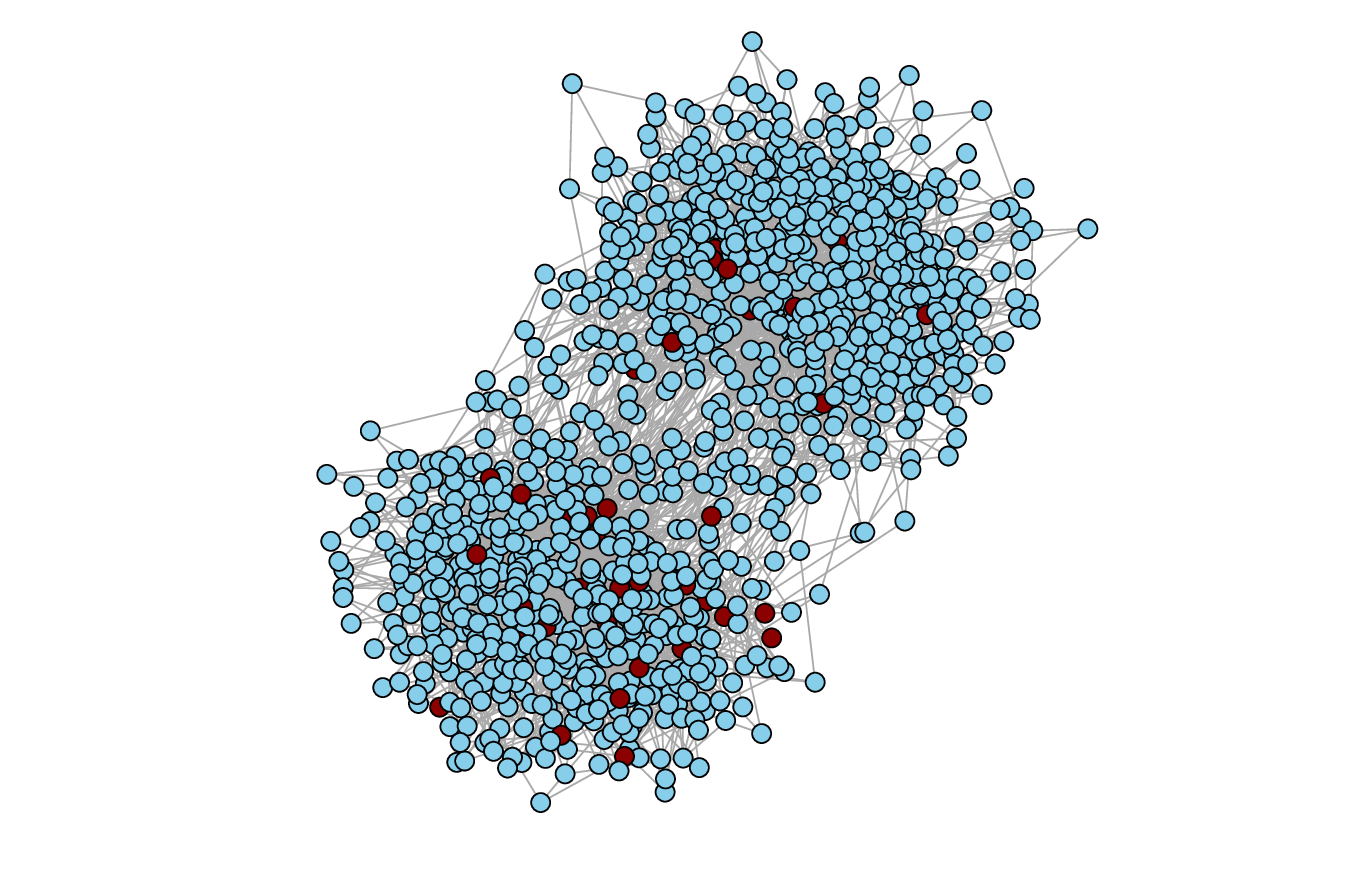}}
	\end{tabular}
	\caption{Examples of simulations with two weakly linked communities and $m=3, p_0=p_1=0,p_2=p_3=1$.  In the left simulation different types are dominating in the two communities, with proportions of red vertices in the two communities being $0.877$ and $0.077$.  In the right simulation blue appears to be dominating in both communities, with proportions of red vertices in the two communities being $0.057$ and $0.021$.\label{commsdiff}}
\end{figure}

Examples of simulations with 1000 vertices and these parameters, with $\theta=0.05$, one with different types dominating in the two communities and one with the same type dominating in both, are shown in Figure \ref{commsdiff}.  The graphs were created using the \texttt{igraph} package in R \cite{igraph}, and plotted using the Fruchterman-Reingold layout, which with these parameters naturally shows the communities.

\subsection{Symmetric community structure, random visible type}\label{randvis}

Here we consider the same community structure as in section \ref{majwins} but the new vertex now picks one of the types it is exposed to uniformly at random.  Hence $p_0=0$, $p_1=p_2=\frac12$ and $p_3=1$, with $R(z)=\frac32(z(1-z)^2+z^2(1-z))+z^3$, with a stable fixed point at $z=\frac12$ and unstable fixed points at $z=0$ and $z=1$, see Figure \ref{fig:4p2}.  In this case the results of \cite{AMR} show that in the single community model the proportion of red tends to $\frac12$ almost surely, and a corresponding result was shown to hold for any number of types and any $m\geq 3$ in \cite{threetypes}.  This rule is increasing, and the only fixed points of $R(z)$ in $[0,1]$ are $z=0$, $z=1/2$ and $z=1$, with the first and last being linearly unstable ($R'(0)=R'(1)=3/2$). Consequently Theorem \ref{same1} shows that the same will apply for the multi community model, regardless of the community structure. Analysis of the stationary points in this case with community structure given by \[A_{\theta}=\begin{pmatrix}1 & \theta \\ \theta & 1\end{pmatrix}\] also appears in \cite{yarrowthesis}, again for $\mu_1=\mu_2=\frac12$, and as might be expected from Theorem \ref{same1} the only stable stationary point is $(\frac14,\frac14,\frac14,\frac14)$, regardless of the value of $\theta$.  So in this case there is no phase transition in $\theta$.

\subsection{Only one possible limit with a single community, different limits in different communities}\label{minority}

\begin{figure}[ht]
	\centering
	\scalebox{0.7}{\includegraphics[trim=6cm 2cm 6cm 1cm, clip=true]{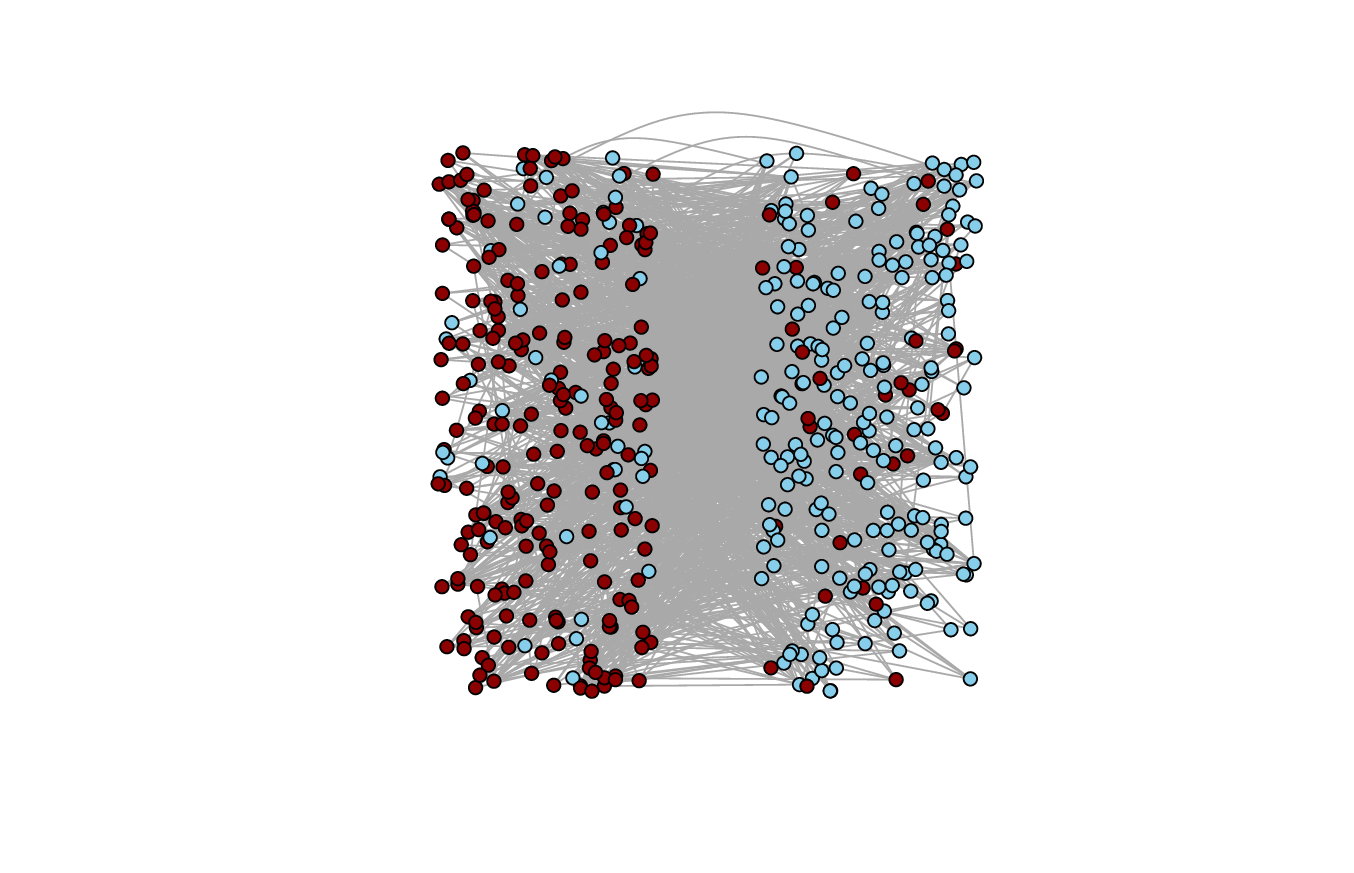}}
	\caption{Example simulation showing different type dynamics in two communities in a case where a single community model converges to a deterministic limit.  The left community has a proportion $0.840$ of its vertices being red, and the right community a proportion $0.151$ of its vertices being red. \label{minoritya10}}
\end{figure}

Theorem \ref{determinant}, part 2, implies that there are examples which that result applies to with the community structure given by \[A=\begin{pmatrix} 1 & \theta \\ \theta & 1\end{pmatrix},\] but with $\theta>1$, indicating that vertices in fact prefer to connect to vertices in the \emph{other} community.  For example, consider this community structure with $\theta=10$, and take $m=3$ with the ``minority rule'': the new vertex chooses to take the type which is less common among its neighbours, giving $p_k=1$ for $k=0,1$ and $p_k=0$ for $k=2,3$.  Here $R(z)=(1-z)^3+3z(1-z)^2$ with a single, stable, fixed point at $z=\frac12$ (see Figure \ref{fig:4p3}), so for the single community model the proportion of red converges to $\frac12$ almost surely.  However, with the above community structure this rule is covered by part 2 of Theorem \ref{determinant} and so shows positive probability of different limits in the two communities.  An example of a simulation with 500 vertices and this rule appears in Figure \ref{minoritya10}, showing different types dominating in the two communities.

\subsection{Example with a touchpoint and an invisible community}\label{tp-rule}

We consider the rule with $m=3, p_0=1/4, p_1=0, p_2=p_3=1$; for this case $R(z)=\frac14(1-z)^3+3z^2(1-z)+z^3$ and has a touchpoint at $z=\frac13$ and a stable fixed point at $z=1$ (see Figure \ref{fig:4p4}), so in the single community model the results of \cite{AMR} will give both $\frac13$ and $1$ as possible limits for the proportion of red vertices.  Proposition \ref{touchpoint} shows that if the community structure matrix $A$ has all entries positive then for small $\theta$ we do not get convergence to different limits here.  However, if one community is invisible to the other this does not apply.

\begin{figure}
	\centering
	\scalebox{0.7}{\includegraphics[trim=6cm 2cm 6cm 2cm, clip=true]{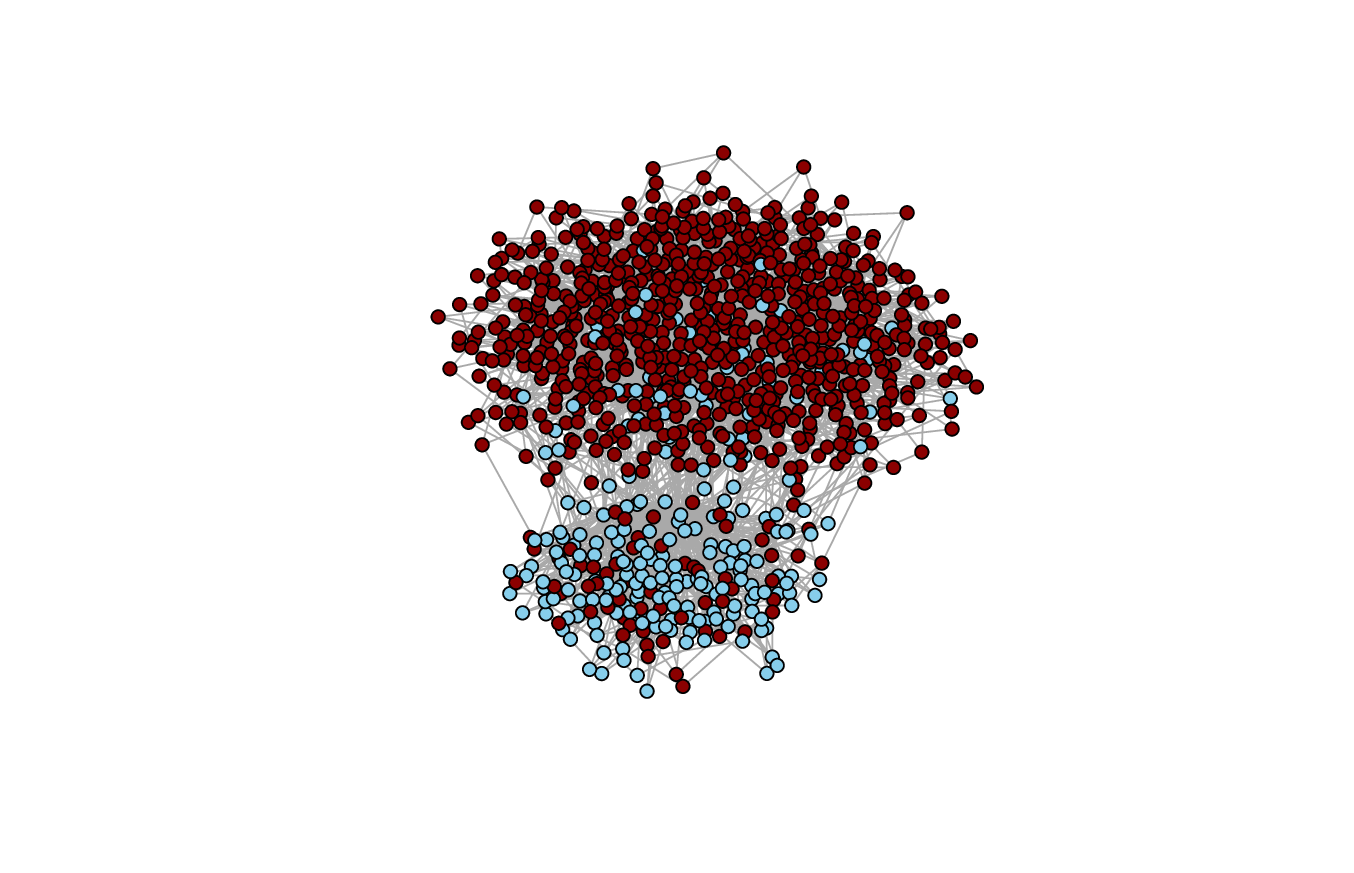}}
	\caption{Here the larger invisible community has a proportion $0.923$ of its vertices red, while the smaller visible community has a proportion $0.269$ of its vertices red, suggesting that the latter may be converging to the touchpoint at $1/3$. \label{touchpoint_invisible}}
\end{figure}

Consider an asymmetric community structure where, in the notation of Theorem \ref{different} \[A_{\theta}=\begin{pmatrix} 1 & 0 \\ \theta & 1\end{pmatrix},\] with $\theta>0$, so that one community (the one labelled 1) is invisible to the other, but the other community is visible to both.  We take $\mu_1=\frac45$ and $\mu_2=\frac15$, so that the invisible community is larger. Figure \ref{touchpoint_invisible} shows an example of a simulation with the type assignment rule described in the last paragraph and this community structure; in this simulation the visible community is close to the touchpoint, but the invisible community is close to the stable root.  Note that it intuitively makes sense that when the invisible community is relatively small it is likely to take the visible community's type, but when it is larger, as here, it is easier for it to maintain its own type dynamics.

\section*{Acknowledgements}
J.H. was supported during this project by the UK Research and Innovation Future Leaders Fellowship MR/S016325/1 and by the European Research Council under the European Union's Horizon 2020 research and innovation programme (grant agreement no.\ 883810). 

Part of this research was carried out during two visits of J.H. to the School of Mathematics and Statistics (SoMaS) in the University of Sheffield, funded by the Heilbronn Institute for Mathematical Research (HIMR) via a Heilbronn Small Grant award. We are grateful to HIMR for this support and to SoMaS for its hospitality.

\end{document}